\documentclass[letterpaper, 10 pt, journal]{IEEEtran}  

\IEEEoverridecommandlockouts                              

\usepackage{lineno,hyperref}
\usepackage{graphicx}
\usepackage{amssymb}
\usepackage{color}
\usepackage{epstopdf}
\usepackage{float}
\usepackage{soul}
\usepackage[utf8]{inputenc}
\usepackage[compress]{cite}
\usepackage{hyperref}
\usepackage{algorithm}
\usepackage{algorithmic,eqparbox}
\usepackage{mathtools}
\usepackage{subcaption}

\usepackage{environ}
\usepackage{tikz}
\usepackage{stackengine}
\usepackage{comment}
\usepackage[normalem]{ulem}

\usepackage{amsmath}         
{
	\newtheorem{theorem}{Theorem}
	
	\newtheorem{assumption}{Assumption}
	\newtheorem{remark}{Remark}

	\newtheorem{lemma}{Lemma}
	
\definecolor{themisGreen}{rgb}{0.0, 0.5, 0.0}
\definecolor{byzantine}{rgb}{0.74, 0.2, 0.64}

\definecolor{themisblue}{rgb}{0.76, 0.89, 1.0}

\newenvironment{list4}{
	\begin{list}{$\bullet$}{%
			\setlength{\itemsep}{0.05cm}
			\setlength{\labelsep}{0.2cm}
			\setlength{\labelwidth}{0.3cm}
			\setlength{\parsep}{0in} 
			\setlength{\parskip}{0in}
			\setlength{\topsep}{0in} 
			\setlength{\partopsep}{0in}
			\setlength{\leftmargin}{0.2in}}}
	{\end{list}}

\title{\huge  {Internal and String Stability of an Observer-based Controller for Vehicle Platooning under the MPF Topology}
}

\author{Wei Jiang, Elham Abolfazli and Themistoklis Charalambous,~\IEEEmembership{Senior Member, IEEE}
	\thanks{W. Jiang and E. Abolfazli are with the Department of Electrical Engineering and Automation, School of Electrical Engineering, Aalto University, 02150, Espoo, Finland (emails: {\tt \{wei.jiang, elham.abolfazli\}@aalto.fi}).}
	\thanks{T. Charalambous is with the Department of Electrical and Computer Engineering, School of Engineering, University of Cyprus, Nicosia, Cyprus. He is also a Visiting Professor at the Department of Electrical Engineering and Automation, School of Electrical Engineering, Aalto University, Espoo, Finland. 
		Email: {\tt charalambous.themistoklis@ucy.ac.cy}.}
	\thanks{The authors have presented a preliminary version~\cite{jiang_platooning} in which only the numeric method by trials and errors is used to design observer-based controller parameters, whereas in this paper two theoretical algorithms {(given the time headway is fixed or for the purpose of minimization)} that are based on a new calculation mechanism are proposed to design controller parameters that guarantee string stability. {Furthermore, how to avoid the peaking phenomenon in the observer-based controller is solved by adopting an additional rule for designing the above parameters.}}
}

\begin{document}

	\maketitle

	%
	%
	%
	%
	\begin{abstract}
		In this paper, we study the internal stability and string stability of a vehicle platoon 
		under the constant time headway spacing (CTHS) policy and the multiple-predecessor-following (MPF) vehicle-to-vehicle information flow topology.  More specifically, we depart from the conventional Proportional-Integral-Derivative (PID) controller design for such systems and we propose the design of an observer-based controller. For designing our observer-based controller, we first design a distributed observer, with which each follower estimates their position, speed and acceleration error with respect to the leader. The observer is designed by means of constructing an observer matrix whose parameters should be determined. Next, we simplify the design of the matrix of the observer in such a way that the design boils down to choosing a single scalar value; this design further simplifies the structure of the controller, whose simplicity enables the derivation of string stability conditions by means of a frequency response method. Subsequently, the string stability conditions for a given time headway, are transformed to conditions for the controller parameters. We obtain controller parameters that satisfy the stability conditions by designing a novel heuristic search algorithm. Furthermore, we extend the search algorithm by incorporating a bisection-like algorithm, which allows to obtain (within some deviation tolerance) the minimum available value of the time headway.
		{Finally, we provide insights about how to finalize the observer-based controller parameters from above algorithms to avoid the peaking phenomenon. }
		The performance of the proposed observer-based controller, which guarantees internal and string stability, is demonstrated via illustrative examples. Additionally, a comparison with a {widely-used} PID controller {for MPF topology} shows that our proposed observer-based controller has better convergence performance. 
		{Also, the platoon safety and controller rate convergence are discussed.}
	\end{abstract}
	

	%
	%
	%
	%
	\section{Introduction}
	
	
	Platoon
	systems in transportation networks refer to  a class of network{ed} systems, in which
	automated vehicles, typically arranged in a string, cooperate via some distributed control protocols, or coupling,  to proceed along the longitudinal direction \cite{levine1966optimal}. Vehicle platooning can boost road throughput and travel safety, while reducing travel time, fuel consumption and \text{CO$_2$} emissions due to the reduced air drag.
	Typically, the distributed control protocol that is designed for a platoon must  guarantee internal stability and  string stability of the platoon system, e.g., see \cite{ploeg2013lp,seiler2004disturbance,middleton2010string,bian2019reducing}.  Internal stability refers to a vehicle matching its speed to the {speed} of the vehicle in front of it, or the platoon leader, while keeping a desired inter-vehicle distance \cite{FENG201981}. String stability enables vehicle platoons  to attenuate the error signals  as they propagate  down the vehicle string \cite{STUDLI2017157}. To be more precise, if the system is string stable, then: (i) vehicles can attain and keep the desired configuration; (ii) the effects of disturbances are attenuated along the string \cite{monteil2019string}. Historically, the work on string stability can be traced back to \cite{peppard1974string} and to the \textit{California PATH} program \cite{sheikholeslam1990longitudinal}.
	
	String stability is known to be dramatically influenced by the spacing policy, which is one of the main components of a platoon that determines the desired inter-vehicle distance. There are two main spacing policies: constant {distance spacing} policy ({CDS}) and constant time headway {spacing} (CTH{S}) policy. Information flow topology (IFT) is another important component  that may impact string stability.
	While a platoon is string unstable with {CDS} policy  under the predecessor-following (PF)  IFT \cite{seiler2004disturbance} or the bidirectional (BD) IFT \cite{barooah2005error}, it is possible to achieve string stability  under the leader-predecessor-following (LPF) \cite{9084094} IFT with the aid of vehicle-to-vehicle (V2V) communication  technology. However, in the LPF topology, each follower needs to know the leader's information via communication channels, which {are often unreliable and become} a  liability, especially as the platoon size grows.
	The CTHS policy, on the other hand, is another way to ensure string stability under the PF topology. In this policy, each vehicle regulates its desired distance from its predecessor by using a linear function of speed (leader/predecessor speed or its own speed) with a constant time headway as the proportional gain \cite{naus2010string}. Apart from the advantages of the CTHS policy, it may compromise the transport throughput, since by using the time headway, the inter-vehicle distance increases as the {speed} grows. It is therefore desirable to reduce the time headway while guaranteeing string stability \cite{flores2018fractional}. The work {in} \cite{8463512} demonstrates that by increasing the number of connected predecessors,  as in the multiple-predecessor-following (MPF) topology{ (an example is shown in Fig.~\ref{fig_platoon})}, the minimum employable time headway will be decreased.
	{Another advantage of adopting MPF IFT is that studies have shown that a car with a velocity of $ 80 \text{km/h} $ following only one predecessor at $ 25 \text{m} $ achieves a $ 30\% $ reduction in aerodynamic drag, and a $ 40\% $ reduction can be attained by following two predecessors~\cite{7286902}.}
	Authors in \cite{bian2019reducing}  consider the CTH{S} policy, the MPF topology and  by introducing a PID controller, they derive the lower bound of the time headway as a function of the number of connected predecessors and the time lag and can also guarantee string stability,.
	{Inspired by the PID controller in~\cite{bian2019reducing} and MPF IFT, authors in~\cite{9301227} designed a similar PID controller with multiple predecessors and  multiple followers topology to deal with limited communication range, time-varying communication delays, and random lossy links.}
	The CTHS policy  and MPF IFT are also adopted in this paper. For a more detailed literature review about the above policies and topologies, please refer to \cite{bian2019reducing,9301227}.
	
	There are mainly {two} analysis methods being used to derive sufficient conditions for proving string stability: $ s $-domain and time-domain analysis methods.  {For nonlinear systems, the time-domain method is usually adopted by using the techniques of Lyapunov functions and eigenvalue analysis; e.g., see \cite{besselink2017string,monteil2019string}. For linear systems, $ s $-domain method{ologie}s  are frequently employed;
		see, e.g., \cite{naus2010string,bian2019reducing,xiao2011practical,abolfazli2020reducing,9462542}.} For a thorough discussion of the relationships and comparisons of the above three methods, the survey paper \cite{FENG201981} is recommended. 
	

	In the literature, all works about vehicle platooning, except for a few, consider  the leader with a constant {speed} {theoretically and mathematically, though they verify their controllers/algorithms for time-varying speed in simulations}. This might lead to reduced performance  when the leader vehicle has speed changes, fluctuations or disturbances, which is how a  vehicle moves in reality. In the literature, the problem of a platoon with a leader whose speed changes over time {theoretically and mathematically} has been studied in two main forms: (i) virtual leader~\cite{monteil2019string,7014426} and (ii) real leader, i.e., leader with specific dynamics~\cite{petrillo2018adaptive,petrillo2020secure,ge2021dynamic,karafyllis2021nonlinear}. Authors in \cite{monteil2019string} designed a control protocol allowing to track a desired (possibly non-constant) reference speed, i.e., the leader does not have specific dynamics and is regarded as a state to be communicated to the followers directly.  In~\cite{7014426}, authors proposed a truck-platoon model in which the speed of the virtual truck is required to be known to all vehicles, including the leader of the platoon. The advantage of this method is that the {changing} speed of the leader is known over time {by the following vehicles in the platoon and, as a result, they can meet the string stability conditions. However,  since the leader's changing speed} is global information, it can limit the {applicability of the method, especially} when there is a large number of vehicles. In~\cite{petrillo2018adaptive,petrillo2020secure,ge2021dynamic}, in which  a leader is employed,  the problem of having the leader with the time-varying speed under the CDS policy is studied.  The distributed adaptive PID controllers in~\cite{petrillo2018adaptive,petrillo2020secure} and PID controller in~\cite{ge2021dynamic} are proposed;  however, no guarantees of string stability  are provided. 
	\cite{karafyllis2021nonlinear} presents a nonlinear controller, however, it is for second-order dynamics.
	
	{There are also different methodologies for studying platoons or connected and automated vehicles (CAV). i) One is optimal control to deal with the noise/disturbance. 
		For example, 
		The optimal state-feedback  Linear-Quadratic-Gaussian (LQG) control  is proposed in~\cite{wang2022optimal} to deal with  time-correlated process noises for the platoon.
		However, there is no string stability guarantee.
		ii) One is sliding mode control (SMC) to deal with uncertainties. Specifically, authors in~\cite{7457311} adopted the distributed adaptive SMC method   to deal with acceleration uncertainties; the same method is used in~\cite{8242685} to deal with  uncertain and time-varying communication topologies.
		iii) Also, the robust control method is utilized in~\cite{gao2016robust} to design $ \mathcal{H}_{\infty} $ controller to simultaneously deal with  vehicle model uncertainties and identical communication delays.
		iv) Another is model predictive control (MPC) which delivers an optimal solution to an objective cost function.
		For instance in order to optimize a vehicle platoon’s fuel consumption,
		authors in~\cite{9801548} combine the switching PID feedback control and distributed economic MPC methods via MPF IFT which helps relax the communication requirements compared with LPF. 
		It is worth noting that the leader of constant speed is still the case in~\cite{9801548} and it is stated clearly in the future work of~\cite{9801548} that it is critical to theoretically analyze the impact of the leader’s time-varying speed. 
		One can see that different methodologies target different factors influencing platoon performance. We use the observer method to propose one observer-based controller to deal with the challenge of the leader vehicle of varying speed mathematically and theoretically, also guarantee the string stability simultaneously.
	}

	In this work, we relax the constraint of the leader's speed being constant and instead assume that it may have a varying transient process with an exponentially converging behavior. Although this kind of leader’s speed is still not the case in reality, due to its behavior and non-zero input{, it is the {first} motivation for building the observer-based controller.}
	Considering a  leader showing the described behavior with specified dynamics (see Eq.~\eqref{leader_dynamics_pva}), the objective of this work is {to} increase the performance of a {vehicle} platoon system, under the directed (see Sec.~\ref{sec:preliminaries} for more details about directed graphs) MPF IFT.  On the way to achieve this objective, first a distributed observer for estimating the  position, {speed} and acceleration error between each vehicle and the leader is proposed. Then,  a  controller  is designed based on  this proposed observer. 
	Although under the MPF topology, the leader shares its information via V2V communications with the first few vehicles, other vehicles do not have access to the leader’s information, unlike other works in the literature (e.g., \cite{besselink2017string,9084094,rodonyi2017adaptive,7014426,monteil2019string,hu2020cooperative})
	where all followers need the (virtual) leader vehicle’s information.
	This is the second motivation, in which we design an observer to estimate and provide additional information for the controller.
	%
	%
	Despite the fact that each vehicle must know the number of its predecessors, which is normally done in advance, the proposed controller is fully distributed and hence can be used for large-scale platooning systems. 
	Furthermore, unlike works~\cite{petrillo2018adaptive,petrillo2020secure,ge2021dynamic}, where internal stability is guaranteed but string stability is not theoretically or mathematically guaranteed, the  controller parameters of this paper are designed  in such a way that  the internal {stability} and {the} string stability of the platoon system are guaranteed.
	The main contributions  are in the following.
	\begin{list4}
		\item 
		{In order to improve platooning control performance under the MPF topology, a distributed observer-based controller is proposed, which can also ensure the internal and string stability of the platoon system}
		One simulated example in Sec.~\ref{exampleA} demonstrates that this distributed observer-based controller achieves a better platooning control performance  compared to the distributed PID controller  of~\cite{bian2019reducing}. The key point in analyzing string stability is to propose a new calculation mechanism for deriving string stability conditions.
		
		\item {Since the $ \mathcal{H}_{\infty} $ norm {of the string stability transfer} function becomes more complicated in our case,  as a result of using this observer-based controller rather than the commonly used PID controller, a novel mechanism is proposed for deriving the string stability conditions. This mechanism utilizes a heuristic searching algorithm {(with observer peaking effect\footnote{{The peak effect means that the trajectories of a closed-loop system (e.g., \eqref{x_tilde2})  significantly deviate from the equilibrium position during the initial phase of the stabilization for some non-zero initial conditions~\cite{polyak2016large}. The large deviation is referred to as an overshoot.}} avoided)} that determines the range of the observer-based controller parameters, given a fixed time headway. This mechanism can also be applied to analyze the string stability via   the $ \mathcal{H}_{\infty} $ norm of PID controller based transfer functions.}
		
		\item A bisection-like algorithm is proposed to obtain the minimum (within a tolerance)  acceptable value of the time headway that can guarantee string stability. Another simulated example in Sec.~\ref{exampleB}  demonstrates that the obtained minimum time headway is smaller than the minimum one from~\cite{bian2019reducing} under the same platooning model conditions.
		
		
	\end{list4}

	The rest of this paper is organized as follows. Sec. \ref{sec:preliminaries} gives some notations and mathematical preliminaries. Sec. \ref{sec:Problem formulation} presents the vehicle model, the spacing policy and the control objective. In Sec. \ref{sec:Main results}, a distributed observer-based controller is proposed with the internal stability   analysis. 
	String stability of the above controller is demonstrated in Sec. \ref{ss}.
	{The controller convergence rate is discussed in Sec.~\ref{sec_convergence_analysis}.}
	Some corroborating simulations are provided in Sec. \ref{sec:examples}. Finally, in Sec. \ref{sec:conclusions} we conclude the paper and discuss future directions.
	
	%
	%
	%
	%
	\section{Notation}\label{sec:preliminaries}
	
	
	$\mathbb{R}^{m \times n}$ and $\mathbb{R}^{n}$ are respectively the $m \times n$ real matrix space and  $ n $-dimensional Euclidean vector space. 
	For the square matrix $ A $,  $ \mathrm{Re}(\lambda(A)) $  represents the real part of eigenvalues of $ A $.
	For any integers $a$ and $b$, with $ a \le b $, denote $ \textbf{I}_a^b=\{a, a+1,\ldots, b\} $.
	The $ \mathcal{H}_{\infty} $ norm of a stable scalar
	transfer function $ H (s) $ is denoted by $ \left\|H\right\|_{\infty} \triangleq \sup_{\omega \in \mathbb{R}} |H(j \omega)| $.
	$ \mathbf{0} $ represents a vector with all elements being 0.

	In a weighted graph $ \mathcal{G} = (\mathcal{N,E,A})$, $\mathcal{N}=\left\{1,2, \ldots, N\right\}$ and $\mathcal{E} \subseteq \mathcal{N \times N}$ are the nodes and edges, respectively. $\mathcal{A} = \left[ a_{ij} \right] \in \mathbb{R} ^{N \times N}$ is the weighted adjacency matrix, where $a_{ij}=1, \left(i,j\right) \in \mathcal{E}$ and $a_{ij} = 0$ otherwise. An edge $\left(j, i\right) \in \mathcal{E}$ means agent $j$ can get information from agent $i$. 
	A directed path from node $i$ to $j$ is a sequence of nodes $ i=l_1, l_2, \ldots, l_t=j $ such that link $ (l_{m+1}, l_{m}) \in \mathcal{E} $ for all $ m=1,2,\ldots, t-1 $.
	The Laplacian matrix $\mathcal{L} = \left[ l_{ij} \right] \in \mathbb{R} ^{N \times N}$ is defined as $l_{ij}= -a_{ij}, i \neq j$ and $l_{ii}= \sum_{j \neq i} a_{ij} $.  
	All nodes that can transmit information to node $i$ directly are said to be in-neighbors of node $i$ and belong to the set $\mathcal{N}^{-}_i=\{ j \in \mathcal{V} \; | \; \varepsilon_{ij} \in \mathcal{E} \}$. 
	The nodes that receive information from node $i$ belong to the set of out-neighbors of node $i$, denoted by $\mathcal{N}^{+}_i=\{ l \in \mathcal{V} \; | \; \varepsilon_{li} \in \mathcal{E} \}$. 
	All in-neighbors and out-neighbors of node $i$ combined are regarded the neighbors of node $i$. 
	
	%
	%
	%
	%
	\section{Problem formulation}\label{sec:Problem formulation}
	\subsection{Longitudinal vehicle platooning dynamics}
	
	The vehicle string is made up of $ N $ follower vehicles with a leader vehicle labeled $ 0 $ and we adopt the following vehicle dynamics~\cite{xiao2011practical,ploeg2013lp,kayacan2017multiobjective}:
	\begin{align}\label{platooning_dynamics}
		\dot{ x}_i(t) = \underbrace{\begin{bmatrix}
				0&1 &0 \\
				0&0 &1 \\
				0&0 &-\frac{1}{\tau}
		\end{bmatrix}}_{\eqqcolon A}   x_i(t) + \underbrace{\begin{bmatrix}
				0\\
				0 \\
				\frac{1}{\tau}
		\end{bmatrix}}_{\eqqcolon B}u_i(t), i \in \textbf{I}_1^N,
	\end{align}
	where $ x_i(t) \coloneqq [ p_i(t),  v_i(t),  a_i(t)]^T, p_i(t), v_i(t), a_i(t) $ represent the longitudinal position, {speed}, and acceleration of vehicle $ i, i\in \textbf{I}_1^N $, respectively;
	$ \tau $ is the engine time constant and $ u_i(t) $ is the vehicle input to be designed;  $ A \in \mathbb{R}^{3\times 3}, B \in \mathbb{R}^{3\times 1} $. 
	One can verify that $ (A, B) $ is controllable.

	The leader's dynamics we investigate here is 
	\begin{equation}\label{leader_dynamics_pva}
		\dot p_0(t) = v_0(t), \,
		\dot v_0(t) = a_0(t),\,
		\dot a_0(t)   =-\frac{1}{\tau}  a_0(t),
	\end{equation}
	which is also
	\begin{equation}\label{leader_dynamics_matrix}
		\dot{ x}_0(t) = A x_0(t).
	\end{equation}
	Note that  many works, e.g.,~\cite{9301227,9801548,bian2019reducing,7134769,ploeg2013graceful}, assume that $ v_0(t) $ is a constant, i.e., $ a_0(t)=0, u_0(t)=0 $, which is quite restrictive. {Authors in~\cite{liu2021cooperative} assume $ \lim_{t \rightarrow \infty} u_0(t) = 0, \lim_{t \rightarrow \infty}  \dot{u}_0(t) = 0 $, which is not easy to relax.}
	In this work, we assume that leader's {speed} has a varying transient, but exponentially converging.
	\subsection{Inter-vehicle distance using constant time headway}
	Since the platooning of homogeneous (same as in~\cite{lunze2018adaptive}) vehicles is investigated here, we set $ h_k =h, \forall k \in \textbf{I}_1^N $ for convenience, where $ h_k $ is the time headway of vehicle $ k $.
	{There are some different desired inter-vehicle distance CTHS policy, e.g., $ d_{i, i-1}(t) = hv_{i-1}(t)+D $,} where $ D $ is the standstill desired gap between vehicle $ i $ and $ i-1 $ (we assume the gap between any two consecutive vehicles are the same since the platoon is homogeneous, to simplify the calculations).
	{As Yanakiev and Kanellakopoulos, in~\cite{728529} stated that ``However, this strategy has a fundamental flaw: if the follower travels at a much higher speed than the leader, say $70$ versus $40~mph$, the desired spacing would be based on the leader’s lower speed; this significantly increases the likelihood of severe collisions'', we adopt the CTHS based on the velocity of the following vehicle ($ v_i(t) $), not the predecessor vehicle ($ v_{i-1}(t) $) as
	}
	\vspace{-0.2cm}
	\begin{equation}\label{distance1}
		d_{i, i-1}(t) = hv_{i}(t)+D.
	\end{equation}
	{This policy \eqref{distance1} is also adopted in many works, e.g.,~\cite{9301227,lunze2018adaptive,728529}}
	One can see \eqref{distance1} also equals $ d_{i, i-j}(t) = \sum_{k=i}^{i-j+1}hv_{k}(t) + jD $, which echos the CTHS policy in~\cite{bian2019reducing}.
	
	\begin{figure}[t]
		\centering
		\includegraphics[width=0.8\columnwidth]{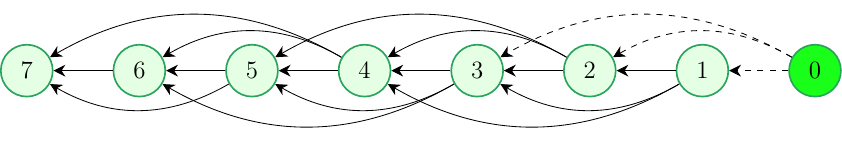}
		\caption{{One example of MPF IFT with seven follower vehicles among which position $ p_i $, velocity $ v_i $, acceleration $ a_i $ and observer  $ \hat{x}_i, i \in \textbf{I}_1^7  $ (see more details in the proposed observer~\eqref{observer_hatx})  are communicated via V2V communication. The leader  sends only its $ p_0, v_0 $ and $ a_0 $ to its three followers. Based on Assumption~\ref{assump_mpf}, we have $ r_1 = 1, r_2 = 2, r_i = r = 3, i \in \textbf{I}_3^7 $.}}
		\label{fig_platoon}
	\end{figure}
	
	\vspace{-0.4cm}
	\subsection{Inter-vehicle communication structure}
	\begin{assumption}\label{assump_mpf}
		The connected vehicles in the platoon are interconnected  via the MPF IFT {with the leader vehicle as the root node}, 
		and the number $ r_i $ of predecessors that follower vehicle $ i $ has {is identical as $ r $, i.e., $ r_i=r\ge 1 $ if $ i\ge r $, and $ r_i = i $ if $ 1 \le i < r, i \in \textbf{I}_1^N  $.}
	\end{assumption}


	{Fig.~\ref{fig_platoon} provides an example showing $ r_i $.}
	Under Assumption~\ref{assump_mpf}, the adjacency matrix $ \mathcal{A} = [a_{ij}] \in  \mathbb{R}^{(N+1)\times (N+1)} $ of  MPF IFT  has the property of $ a_{0j}=0,  a_{ij}=0, \forall i <j, i \in \textbf{I}_1^N, j  \in \textbf{I}_1^N $, i.e., $ \mathcal{A} $ becomes a lower-triangular matrix.
	As a result, the Laplacian matrix $ \mathcal{L} $ has the same property. Partition $ \mathcal{L} $ as 
	$
	\mathcal{L} =  
	\begin{bmatrix}
		0 & 0_{1 \times N} \\
		\mathcal{L}_{N \times 1} & \mathcal{L}_{1} 
	\end{bmatrix} ,
	$ where
	$ \mathcal{L}_1\in \mathbb{R}^{N\times N} $ with  $ l_{ii} = r_i, i \in \textbf{I}_1^N $.

	\vspace{-0.3cm}
	\subsection{Control objectives}\label{obj}
	Define the desired position for the $ i $th  vehicle related to its predecessor  vehicle  as
	\vspace{-0.2cm}
	\begin{align}
		\bar p_i^d(t) = p_{i-1}(t)- d_{i, i-1}(t),
	\end{align}
	Then, define the predecessor-follower position, {speed} and acceleration errors for  the $ i $th  vehicle related to its predecessor  vehicle as follows:
	\vspace{-0.2cm}
	\begin{align}
		\bar{p}_i(t) =&  p_i(t) - \bar{p}_i^d(t) = p_i (t)- p_{i-1}(t) + h{v_{i}}(t)+D, \nonumber\\
		\bar{v}_i(t) =&  v_i(t) - v_{i-1}(t), \nonumber\\
		\bar{a}_i(t) =&  a_i(t) - a_{i-1}(t), i \in \textbf{I}_1^N . \label{pva_error2}
	\end{align}
	Denote an augmented variable $ \bar x_i(t) \coloneqq [\bar p_i(t), \bar v_i(t), \bar a_i(t)]^T $ and under Assumption~\ref{assump_mpf}, the platooning is required to track a  varying-{speed} leader, where the following objectives are set:
	\begin{itemize}
		\item[$\mathrm{O_1}$:] convergence of the predecessor-follower platooning tracking error $ \bar x_i(t) $;
		\item[$\mathrm{O_2}$:] guarantee of the predecessor-follower string stability.
	\end{itemize}

	{In the following, for the convenience of presentation, the time index $ (t) $ is omitted. The distributed observer-based controller is proposed in Sec.~\ref{sec:Main results} to achieve objective $\mathrm{O_1}$. The objective $\mathrm{O_2}$ is achieved in Sec.~\ref{ss}.}

	\vspace{-0.2cm}
	\section{Distributed observer-based controller}\label{sec:Main results}
	In this section, we first provide a leader-following platooning tracking error model which is related to the leader; then, a distributed observer is proposed to estimate this error; thereafter, the link between this error with the predecessor-follower platooning tracking error $ \bar x_i $ in~\eqref{pva_error2} is presented and the convergence is proved to achieve objective $\mathrm{O_1}$.
	\vspace{-0.2cm}
	\subsection{Leader-following platooning tracking error model}
	Define the desired position for the $ i-$th follower vehicle from the leader vehicle as
	\begin{align*}
		p_i^d 
		= p_0 - {i}\left(hv_0 + D\right).
	\end{align*}
	Then, the leader-following position, {speed} and acceleration errors are respectively defined as follows:
	\begin{align}
		\tilde{p}_i =&  p_i - p_i^d = p_i -p_0 + {i}(hv_0 + D), \nonumber\\
		\tilde{v}_i =&  v_i - v_0, \nonumber\\
		\tilde{a}_i =&  a_i - a_0. \label{pva_error}
	\end{align}
	Now, we calculate $ \dot{\tilde{p}} _i  $ first as
	\begin{align}
		\dot{\tilde{p}}  _i  = \tilde{v}_i +  iha_0.\label{Omega}  
	\end{align}
	By defining the augmented variable $ \tilde x_i$ as 
	$ 
	\tilde x_i \coloneqq [\tilde p_i, \tilde v_i, \tilde a_i]^T ,
	$
	from \eqref{platooning_dynamics} we obtain
	\begin{align}\label{platooning_error_dynamics}
		\dot{\tilde x}_i = A  \tilde x_i + B u_i + B_1\Omega_i, i \in \textbf{I}_1^N,
	\end{align}
	where $ B_1 = \begin{bmatrix}
		1&
		0 &
		0
	\end{bmatrix}^T, \Omega_i = iha_0 $.
	Also, from \eqref{leader_dynamics_pva} we have that $ \dot a_0 =- \frac{1}{\tau} a_0 \Rightarrow \lim_{t\rightarrow \infty} a_0(t) =0 $. As a result, we deduce that $\lim_{t\rightarrow \infty} \Omega_i(t) = 0$.
	\subsection{Distributed observer design} 
	Analyzing the construction of the predecessor-follower platooning tracking error $ \bar x_i $ in~\eqref{pva_error2} and the leader-following platooning tracking error $ \tilde x_i $ in~\eqref{pva_error}, we find out that as $ \tilde x_i \to 0, t\rightarrow \infty$, then eventually $v_i=v_j= v_0$ and $a_i=a_j=a_0$. As a result, eventually $\bar{v}_i=0$, $\bar{a}_i =0$. At the same time, $ \tilde x_i \rightarrow 0 \Rightarrow \tilde{p}_i \rightarrow 0,  \tilde{p}_j \rightarrow 0 \Rightarrow  \tilde{p}_i - \tilde{p}_{i-1} = p_i -p_{i-1}+ hv_0+D \rightarrow 0 \Rightarrow \bar{p}_i=0, t\rightarrow \infty $ as $v_i \rightarrow v_0$, i.e., $ \tilde x_i\rightarrow 0\Rightarrow  \bar x_i\rightarrow 0$. Therefore, objective $\mathrm{O_1}$ essentially involves proving the convergence of $ \tilde x_i $.
	In fact, the relation between $ \bar x_i $ and $ \tilde x_i $ is
	\begin{equation}\label{relation_bar_tilde_x}
		\bar x_i = \tilde x_i - \tilde x_{i-1} + \begin{bmatrix}
			h\tilde v_{i}&
			0 & 
			0
		\end{bmatrix}^T, i \in \textbf{I}_1^N, \tilde x_{0}= \mathbf{0}.
	\end{equation}

	The idea is to design a distributed observer as $ \hat{x}_{i} \coloneqq [\hat  p_i, \hat v_i, \hat a_i]^T \in  \mathbb{R}^{3} $ with  $ \hat{x}_{i}(0) =\mathbf{0} $ to estimate the leader-following platooning  error $ \tilde{x}_{i} $. 
	
	From Assumption~\ref{assump_mpf}, the neighbor of  vehicle $1$ is only the leader, which means vehicle $1$ can receive the information of $ p_0, v_0, a_0 $.
	As a consequence, the observer mathematical format is divided into $ \hat{x}_{1} $ and $ \hat{x}_{i},i \in \textbf{I}_2^N $, respectively, {as proposed in~\eqref{observer_hatx},}
			\begin{figure*}[t]
				\begin{subequations}
					\label{observer_hatx}
					\begin{align}
						\dot{\hat{x}}_{1}=&A \hat{x} _{1} + B u_1  
						+
						BK({\begin{bmatrix}
								p_1-p_0+hv_0 +D\\
								v_1-v_0 \\
								a_1-a_0
						\end{bmatrix}} - \hat{x}_1)+ BL  ({\begin{bmatrix}
								p_1-p_0+hv_0 +D\\
								v_1-v_0 \\
								a_1-a_0
						\end{bmatrix}} - \hat{x}_1), \label{observer_hatx1}\\
						\dot{\hat{x}}_{i}=&A \hat{x} _{i} + B u_i +
						BK({\begin{bmatrix}
								p_i-p_{i-1}+hv_{i-1} +D\\
								v_i-v_{i-1} \\
								a_i-a_{i-1}
						\end{bmatrix}} - \hat{x}_i)  + BL  \{\sum_{j=1}^{i-1} a_{ij} [\begin{bmatrix}
							p_i-p_j+{(i-j )} hv_j +(i-j )D\\
							v_i-v_j \\
							a_i-a_j
						\end{bmatrix}- (\hat{x}_i-\hat{x}_j)]  \nonumber \\
						& +a_{i0} ({\begin{bmatrix}
								p_i-p_{i-1}+hv_{i-1} +D\\
								v_i-v_{i-1} \\
								a_i-a_{i-1}
						\end{bmatrix}} - \hat{x}_i)\}, i \in \textbf{I}_2^N, \label{observer_hatx2}
					\end{align}
				\end{subequations}
			\end{figure*}
			where 
			$ a_{ij} $ is the element of the adjacency matrix $ \mathcal{A} $ with $ a_{i0}=1, i \in \textbf{I}_1^r  $ and $ a_{i0}=0, i \in \textbf{I}_{r+1}^N $ {($ r $ is from Assumption~\ref{assump_mpf})}; observer parameters  $ L \in  \mathbb{R}^{1\times 3 }$ and $K \in  \mathbb{R}^{1\times 3 } $ will be designed later.
			Note that $ a_{i0}=0, i \in \textbf{I}_{r+1}^N $ means vehicles $ r+1, r+2, \ldots, N $ do not need to know the leader's information, which is a crucial departure from the controllers in~\cite{9084094}/\cite{7014426} in which the velocity of the leading virtual truck should be known to all vehicles in the platoon.

			One can see that the observer of vehicle $ i $ requires the relative position, the relative acceleration with respect to its neighbor $ j $, and additionally the velocity $ j $ and observer information  $ \hat x_j $ from its neighbor $v_j\in\mathcal{N}_i$. 
			To avoid the excessive communication, a special design of the parameter matrix $ L $ will be presented in Sec.~\ref{sec_L}, in which the information needed will be much less, as it will be demonstrated in Sec.~\ref{sec_ss_condition}.


				
				Denote the observer estimating error $ \xi_i \in  \mathbb{R}^{3} $ as  
				\begin{equation}\label{ESO_error_dynamics}
					\xi_i \coloneqq \tilde{x}_i - \hat{x} _{i} , i \in \textbf{I}_1^N.
				\end{equation}
				\begin{remark}\label{remark_vehicle1}
					Recall that the purpose of designing the observer is to have $ \lim_{t\rightarrow \infty} \xi_i(t) =0 $. 
					From~\eqref{relation_bar_tilde_x},
					one can see $ \bar x_1 = \tilde x_1 $; then, one can design $ \hat x_1 = \tilde x_1 $ directly to have  $ \xi_1 =0 $ for all time. However, with this design, we cannot guarantee the platoon string stability theoretically and mathematically as we will demonstrate this point in Remark~\ref{remark_x1}.
				\end{remark}
				
				After some algebraic manipulations, we obtain that
				\begin{align*}
					p_i-p_j+& (i-j)hv_j +(i-j )D 
					= p_i-p_0+ ihv_0 - [p_j-p_0  \\
					&+jhv_0 + (i-j)hv_0 ]  +(i-j)hv_j +(i-j )D \\
					&= \tilde{p}_i-\tilde{p}_j + (i-j)h\tilde{v}_j.
				\end{align*}
				Hence, the proposed observer $ \hat{x}_{i} $, $ i \in \textbf{I}_1^N$ in \eqref{observer_hatx} changes to 
				\begin{align}
					&\dot{\hat{x}}_{1}= A \hat{x} _{1} + B u_1  
					+
					BK\xi_1+ BL   \xi_1, \nonumber\\
					&\dot{\hat{x}}_{i}=A \hat{x} _{i} + B u_i +
					BK\xi_i + BL  \sum_{j=1}^{i} l_{ij} \xi_j     \label{modified_observer4}\\
					&+ \underbrace{BL \sum_{j=1}^{i-1} a_{ij} \begin{bmatrix}
							(i-j)h\tilde v_{j}\\
							0 \\
							0
						\end{bmatrix}
						-B(a_{i0}L +K)\begin{bmatrix}
							\tilde p_{i-1}-h\tilde v_{i-1}\\
							\tilde v_{i-1} \\
							\tilde a_{i-1}
					\end{bmatrix}}_{\eqqcolon\Pi_{i-1}}. \nonumber
				\end{align}

				Now, based on the leader-following platooning  error dynamics $ \tilde{x}_{i} $ in~\eqref{platooning_error_dynamics},  we design the control input as
				\begin{equation}\label{input}
					u_i= -K \hat{x}_{i}, i \in \textbf{I}_1^N,
				\end{equation}
				such that
				\begin{equation}\label{x_tilde2}
					\dot{\tilde{x}}_{i}  = (A-BK)\tilde{x}_{i}+BK\xi_i+ B_1\Omega_i.
				\end{equation}
				Due to the fact that $ \lim_{t\rightarrow \infty} \Omega_i(t) = 0 $ in~\eqref{Omega}, if $ \lim_{t\rightarrow \infty} \xi_i(t) = 0 $  and we  design $ K $ such that $ A-BK $ is Hurwitz, then, $ \lim_{t\rightarrow \infty} \tilde{x}_{i}(t) = 0 $ {such that $ \lim_{t\rightarrow \infty} \bar{x}_{i}(t) = 0 $}.

				\begin{lemma}\label{theorem_convergece_platooning}
					Under Assumption~\ref{assump_mpf} with $ r_i $ being the number of predecessors of follower vehicle $ i $, the predecessor-follower platooning tracking error will converge to zero asymptotically, i.e., $  \lim_{t\rightarrow \infty} \bar{x}_{i}(t) = 0 $,
					by designing parameter matrices $ K $ and $ L $ such that $ A-BK $, $ A-BK-{r_i} BL $ are Hurwitz.
				\end{lemma}
				\begin{proof}
					See Appendix~\ref{proof:lemma}.
				\end{proof}
				
				\begin{remark}\label{remark_oneAfterOne}
					The proof of Lemma~\ref{theorem_convergece_platooning} states that the platooning tracking error of vehicle $ 1 $ will converge first, then comes the convergence of vehicle $ 2 $, then vehicle $ 3 $, \ldots, and finally the convergence of vehicle $ N $. This is reasonable as in vehicle platooning, each vehicle $ i $ has at least its predecessor vehicle $ i-1 $ as its neighbor and needs the information of vehicle $ i-1 $ for observer design, as one can see $ p_{i-1}, v_{i-1}, a_{i-1} $ inside the proposed observer $  \hat{x} _{i} $~\eqref{observer_hatx}.
				\end{remark}
				
				In Lemma~\ref{theorem_convergece_platooning}, for the observer parameter matrices design, it is trivial to design $ K $ for $ A-BK $ being Hurwitz. In addition to that, how to design $ L $ to have $ A-BK- {r_i} BL $ be Hurwitz deserves special attention.
				\begin{remark}
					We do not design $ L= \mathbf{0}_{1\times 3 } $ such that 	$ A-BK- r_i BL $ becomes $ A-BK $ to be Hurwitz. One reason is that based on our proposed controller $ u_i $ in~\eqref{input} and observer $ \hat{x}_i $ in~\eqref{observer_hatx},  if $ L =  \mathbf{0}_{1\times 3 } $, then, our observer based controller becomes a controller for PF IFT as the information exchange term among multiple neighbors will disappear in~\eqref{observer_hatx}. The other reason is that based on \cite[Theorem 2]{bian2019reducing} and \cite[Theorem 2]{9462542}, the existence of multiple predecessors in MPF IFT can decrease the lower bound of time headway, i.e., the larger the number of predecessors, the smaller the value of the lower bound of time headway. By  designing $ L \neq \mathbf{0}_{1\times 3 } $, it means our controller is proposed for MPF IFT aiming at decreasing the lower bound of time headway.
				\end{remark}
				
				\vspace{-0.4cm}
				\subsection{Parameter matrix $ L $ design}\label{sec_L}

				%
				
				Since $ A-BK $ is Hurwitz, i.e., $ \mathrm{Re}(\lambda(A-BK) )<0 $ and $ {r_i} $ is positive from Assumption~\ref{assump_mpf}, inspired by~\cite{9130792}, 
				we design  the term $ BL $ being non-negative definite by proposing one solution of $ L $ as
				\vspace{-0.3cm}
				\begin{equation}\label{L}
					L = \alpha B^T, \alpha >0.
				\end{equation}
				Thus, the design of matrix $ L $ is simplified to design the scalar $ \alpha $, which will be illustrated in Sec.~\ref{sec_ss_condition}.
				
								
								\subsection{Platoon safety}
								
								{As the extensive study of string stability lacks the safety analysis \cite[Sec. I]{7547317}, we try to include it in this paper by introducing some surrogate safety measures (SSMs).
									According to the latest survey \cite[Fig. 2]{10052670}, time-based SSMs are the most frequently used SSM for mixed traffic safety assessment, followed by deceleration-based and distance-based SSMs. 
									There is no generally accepted consensus to classify acceptable and non-acceptable risk level produced by different SSMs \cite[Sec. 1]{LU2021106403}.
									Authors in~\cite{10052670} also state that no collective
									guidance exists to propose the best set of SSMs that can accurately report traffic conflicts. Consequently,
									inspired by the simulation in \cite[Fig. 4 (c)]{10052670}, we choose three SSMs, i.e., time to collision (TTC), deceleration rate to avoid crash (DRAC) and difference of the space distance and stopping distance (DSS) accordingly. }
								
								{In detail, TTC and DRAC are defined in \cite[Equations (1) and (4)]{LU2021106403}, respectively;
									DSS is defined in \cite[Equation (4)]{okamura2011impact}.
									These three SSMs are practical and easy to calculate. Denote their thresholds respectively as $ TTC^{*}, DRAC^{*} $ and 0.
									It is worth noting that
									no consensus exists for selecting the SSM thresholds; nonetheless, the outcome of safety analysis relies heavily on threshold selection~\cite{10052670}.
									Although string stability strictly speaking is not considered as SSM, previous literature has proved that better string stability could bring important safety benefits \cite[Sec. 3.3]{WANG2021106157}. As a result, since string stability is considered and solved in this paper (see Sec.~\ref{ss} for more details), we choose less conservative thresholds, e.g., $ TTC^{*} = 2s, DRAC^{*}= 3.4 m/s^2 $ ($  TTC^{*} $ range is typically $ 1.5s-2s $; $ DRAC^{*} $ range is $ 3m/s^2-3.4 m/s^2 $~\cite[Sec.~III-C]{10052670}).}
								
								{When the following three unsafe conditions:
									\begin{equation}\label{safety_condition}
										TTC \le TTC^{*}, DRAC\ge DRAC^{*}, DSS \le 0
									\end{equation}
									are satisfied simultaneously, we mark the situation as unsafe under our observer-based controller.
									Consequently, the initial condition of the platoon system should not satisfy the unsafe condition~\eqref{safety_condition}.}
								{Note that the underline relationship between string stability and individual vehicle's safety is not explicitly understood and needs to be further explored~\cite[Sec. 4.4]{WANG2021106157}, which is one of our future research directions.}



								\vspace{-0.3cm}
								\section{String stability}\label{ss}
								
								In this section, we first present the string stability condition related to the observer-based controller parameters. Then, a new mechanism is proposed to determine how to design these parameters from the {string stability} condition using a heuristic searching algorithm. Finally, a time headway minimization algorithm is proposed, which  demonstrates that the proposed observer-based controller is capable of further improving the traffic throughput by simulations in~Sec.~\ref{exampleB}.
								
								\vspace{-0.4cm}
								\subsection{Conditions for string stability}\label{sec_ss_condition}
								
								Typically, the variation of the leading vehicle’s {speed} is viewed as a disturbance on the platoon, which results in a certain transient process. The property of this transient process is studied by using the notion of string stability.
								
								For the string stability,
								we consider the amplification of spacing errors since spacing errors directly affect the platoon safety. Define the predecessor-follower spacing error as 
								\begin{equation}\label{bar_e_i}
									\bar e_i\coloneqq \bar{p}_i =p_i - p_{i-1} + h v_{i} +D ,
								\end{equation}
								where $ \bar{p}_i $ is defined in~\eqref{pva_error2}. In order to prove string stability we require that: $\left\| \bar e_i \right\|_{\infty} \le \left\| \bar e_{i-1} \right\|_{\infty} $.
								
								
								
								Based on Assumption~\ref{assump_mpf},  observer $\hat{x} _{i}(t),  i \in \textbf{I}_1^N $ in \eqref{observer_hatx} can be transformed into~\eqref{observer_r}.
								\begin{figure*}[t]
									\begin{align}
										\dot{\hat{x}}_{i}=&A \hat{x} _{i}+ B u_i +
										BK({\begin{bmatrix}
												p_i-p_{i-1}+hv_{i-1} +D\\
												v_i-v_{i-1} \\
												a_i-a_{i-1}
										\end{bmatrix}} - \hat{x}_i)  + BL  \sum_{l=1}^{r} [\begin{bmatrix}
											p_i-p_{i-l}+lhv_{i-l} + lD\\
											v_i-v_{i-l} \\
											a_i-a_{i-l}
										\end{bmatrix} - (\hat{x}_i-\hat{x}_{i-l})]. \label{observer_r}
									\end{align}
								\end{figure*}
								Here, define $ \hat{x}_0=0 $ in case $ l=i $.
								
								From~\eqref{L}, we know $ L = \alpha B^T $. 
										Denote
										\begin{equation}\label{design_k}
											K = \begin{bmatrix}
												k_1&k_2 &k_3 
											\end{bmatrix},
										\end{equation}
										where scalars $ k_1, k_2, k_3  $ are to be decided. Then, based on input~\eqref{input}, we calculate the matrix form of $ BL, BK, A-BK $ such that \eqref{observer_r} can change to 
											\begin{align}
												\begin{bmatrix}
													\dot{\hat p}_i\\
													\dot{\hat v}_i \\
													\dot{\hat a}_i
												\end{bmatrix} = & \begin{bmatrix}
													0&1 &0 \\
													0&0 &1 \\
													-\frac{k_1}{\tau}&-\frac{k_2}{\tau} &-\frac{1+k_3}{\tau}
												\end{bmatrix}  \begin{bmatrix}
													\hat p_i\\
													\hat v_i \\
													\hat a_i
												\end{bmatrix} +   \sum_{l=1}^{r}\begin{bmatrix}
													0&0 &0 \\
													0&0 &0 \\
													0&0 &\frac{\alpha}{\tau^{2}}
												\end{bmatrix} \nonumber\\ & \times (\begin{bmatrix}
													p_i-p_{i-l}+lhv_{i-l} + lD\\
													v_i-v_{i-l} \\
													a_i-a_{i-l}
												\end{bmatrix} - \begin{bmatrix}
													\hat p_i- \hat p_{i-l}\\
													\hat v_i- \hat v_{i-l} \\
													\hat a_i- \hat a_{i-l}
												\end{bmatrix})
												\nonumber\\
												&+
												\begin{bmatrix}
													0&0 &0 \\
													0&0 &0 \\
													\frac{k_1}{\tau}&\frac{k_2}{\tau} &\frac{k_3}{\tau}
												\end{bmatrix} \begin{bmatrix}
													p_i-p_{i-1}+hv_{i-1}+D-\hat p_i\\
													v_i-v_{i-1}-\hat v_i \\
													a_i-a_{i-1}-\hat a_i
												\end{bmatrix},
											\end{align}
											which is also equivalently written as
											\begin{align}
												\dot{\hat p}_i=&\hat{v}_i, 
												\dot{\hat v}_i=\hat{a}_i, \nonumber\\
												\dot{\hat a}_i =& -\frac{k_1}{\tau}\hat p_i-\frac{k_2}{\tau}\hat v_i-\frac{1+k_3}{\tau}\hat a_i  
												+\frac{k_1}{\tau}(p_i-p_{i-1}+hv_{i-1}\nonumber\\ &+D -\hat p_i)+\frac{k_2}{\tau} (v_i - v_{i-1} -\hat v_i)+\frac{k_3}{\tau}(a_i - a_{i-1} -\hat a_i) \nonumber\\&
												+  \sum_{l=1}^{r}\frac{\alpha}{\tau^{2}}[a_i-a_{i-l}-(\hat a_i- \hat a_{i-l})].
												\label{hat_a_i}
											\end{align}
											\begin{remark}
												%
												One can see 
												observer $\hat{x} _{i}  $ in~\eqref{observer_hatx} actually has a third order integrator dynamics. After $ L $ is designed in~\eqref{L}, the resulted observer only needs the relative position with respect to its predecessor, velocity of its predecessor, relative accelerations with respect to its neighbors and acceleration observer of neighbors, which accounts for less {communicated} information compared to the original observer~\eqref{observer_hatx}.
											\end{remark}
											\begin{theorem}\label{inter_string_theorem}
												{Under Assumptions~\ref{assump_mpf}}, the internal and string stability of the platooning system~\eqref{platooning_dynamics} and \eqref{leader_dynamics_pva} is guaranteed by the observer-based controller \eqref{hat_a_i} and \eqref{input} if the controller parameters $ k_1, k_2, k_3, \alpha $ satisfy $ A-BK $ is  Hurwitz with $ K = \begin{bmatrix}
													k_1&k_2 &k_3 
												\end{bmatrix} $  and
												\begin{equation}\label{string_stability_condition}
													\left\|H(j \omega)\right\|_{\infty}=
													\left\| \frac{q_1(j\omega) T_4}{T_1T_3+ T_2T_4} \right\|_{\infty} \leq 1, \forall \omega \in [0, \infty),
												\end{equation}
												where 
												\begin{subequations}
													\label{T1234}
													\begin{align}
														T_1 =& 2k_1+2k_2s + (1+2k_3+ r\bar \alpha )s^2+\tau s^3, \\
														T_2 =& k_1+ k_2s + (k_3+ r\bar \alpha )s^2,\\
														T_3 =& \tau s^3+ s^2,\\
														T_4 =& k_1+k_2s+k_3s^2,\\
														q_1(s) =& \bar{\alpha} s^2 +k_1 - (k_1h-k_2) s+k_3s^2, \label{q1_equation}\\
														\bar{\alpha}=& \frac{\alpha}{\tau}, \alpha >0.
													\end{align}
												\end{subequations}
											\end{theorem}
											\begin{proof}
												See Appendix~\ref{proof:theorem:inter_string}.
											\end{proof}
											
											\vspace{-0.2cm}
											\subsection{Parameters $ k_1, k_2, k_3 $ and $ \alpha $ for string stability}
											{In the literature, only a few works consider a platoon with a leader whose speed changes over time theoretically and mathematically, e.g., ~\cite{petrillo2018adaptive,petrillo2020secure,ge2021dynamic}. However, no guarantees of string stability  are provided in the above works.}
											{In this subsection, we propose a new string stability condition parameter designing mechanism.}

											One can see that string stability condition \eqref{string_stability_condition} is quite complicated and it is not obvious how to design $ \alpha $ and $ k_1,k_2,k_3 $. Further analysis of the structure of $ \left\|H(j \omega)\right\|_{\infty} $ is needed. Towards this end, from \eqref{T1234} we deduce that 
											\begin{align*}
												q_1(s) =& \bar \alpha s^2 - k_1hs + T_4, \\
												T_1 =& \tau s^3 + (1+r\bar \alpha) s^2 + 2T_4, \\
												T_2 =& r\bar \alpha s^2 + T_4.
											\end{align*}
											Thus, the numerator and denominator of $ H(s) $ change respectively to
											\begin{align*}
												q_1(s)T_4 =& s(\bar \alpha s - k_1h)T_4 + T_4^2, \\
												T_1T_3+ T_2T_4 
												=&  s^4(\tau s +1) (\tau s +1 +r\bar \alpha)  + s[2\tau s^2 \nonumber\\ &+ (2 + (r-1)\bar \alpha)s +  k_1h] T_4 + q_1(s)T_4.  \nonumber
											\end{align*}
											After some algebraic manipulation, {we propose a new structure of} $ H(s) $ as
											\begin{align}
												H(s) =& \frac{1}{X + 1},\label{eq_X_H}\\
												X=& \frac{T_1T_3+ T_2T_4 -  q_1(s)T_4}{ q_1(s)T_4}.\label{eq_X}
											\end{align}
											Since $ X $ is a complex number, so it can be written as $ X \coloneqq \mathrm{Re}(X) + j\mathrm{Im}(X) $. Therefore, if $ \mathrm{Re}(X) \ge 0 $ or $ \mathrm{Re}(X) \le -2 $, then we have $ \mathrm{Re}(X+1) \ge 1 $ or $ \mathrm{Re}(X+1) \le -1 $, respectively. As a consequence, whatever $ \mathrm{Im}(X+1) $ is, $ \left\|X+1\right\|_{\infty} \ge 1 \Rightarrow \left\|H(j\omega)\right\|_{\infty} \le 1 $. 
											
											To sum up here, by the above transformations, our focus for the string stability switches from \eqref{string_stability_condition} to \eqref{eq_X}.
											When $ s = j\omega $, denote
											\begin{align}
												X &\coloneqq  \frac{\mathrm{Re}(X_{\text{num}}) + j\mathrm{Im}(X_{\text{num}})}{\mathrm{Re}(X_{\text{den}}) + j\mathrm{Im}(X_{\text{den}}) }\\ & = \frac{ (\mathrm{Re}(X_{\text{num}}) + j\mathrm{Im}(X_{\text{num}}) )( \mathrm{Re}(X_{\text{den}}) -  j\mathrm{Im}(X_{\text{den}}) ) }{\underbrace{\mathrm{Re}(X_{\text{den}})^2 + \mathrm{Im}(X_{\text{den}}) ^2}_{\eqqcolon Z}} \nonumber \\
												&=  \underbrace{(\mathrm{Re}(X_{\text{num}}) \mathrm{Re}(X_{\text{den}}) + \mathrm{Im}(X_{\text{num}})   \mathrm{Im}(X_{\text{den}}))}_{\eqqcolon Y}  /Z  + j\mathrm{Im}(X). \nonumber
											\end{align}
											One can see that the value of denominator of $ \mathrm{Re}(X) $ is positive but difficult to be calculated, i.e., to derive $ \mathrm{Re}(X) \le -2 $ is difficult. However, to have $ \mathrm{Re}(X) \ge 0 $, we just need to have the numerator of $ \mathrm{Re}(X) $ to be non-negative, i.e., to have
											\begin{subequations}
												\label{check}
												\begin{align}
													Y \ge& 0, \label{parameter_setting_condition}\\
													\mathrm{Re}(X_{\text{num}}) =& n_2 \omega^2 + n_4 \omega^4 +n_6 \omega^6,\\
													\mathrm{Im}(X_{\text{num}}) =& n_1 \omega^1 + n_3 \omega^3 +n_5 \omega^5,\\
													\mathrm{Re}(X_{\text{den}}) =& d_0 + d_2 \omega^2 + d_4 \omega^4,\\
													\mathrm{Im}(X_{\text{den}}) =& d_1 \omega^1 + d_3 \omega^3,
												\end{align}
											\end{subequations}
											where 
											{\begin{subequations}
													\label{parameters}
													\begin{align}
														n_1 =& hk_1^2, \\
														n_2 =&-2k_1- hk_1k_2 - \bar \alpha k_1(r-1), \\
														n_3 =&  -2k_1\tau -2k_2 - hk_1k_3 - \bar \alpha k_2(r-1),\\
														n_4 =& 2k_2\tau + 2k_3+1  + \bar \alpha r + \bar \alpha k_3(r-1), \\
														n_5 =& 2\tau + \alpha r+ 2k_3\tau, \\ 
														n_6 =& -\tau^2, \label{n_6}
													\end{align}
													and
													\begin{align}
														d_0 =& k_1^2, \\ 
														d_1 =& k_1 (2 k_2-hk_1 ) ,  \label{d_1} \\ 
														d_2 =& -k_2^2 -2 k_1k_3 + hk_1k_2 - \bar \alpha k_1, \label{d_2}\\
														d_3 =& -k_3 (2 k_2-hk_1 ) - \bar \alpha k_2,  \label{d_3} \\ 
														d_4 =& k_3^2 + \bar \alpha k_3.  \label{d_4}
													\end{align}
											\end{subequations}}
											
											\begin{remark}\label{reX_2}
												{From~\eqref{eq_X_H} and \eqref{eq_X}, to have $ \mathrm{Re}(X) \ge 0 $, we just need $ Y\ge 0 $ as in~\eqref{parameter_setting_condition}; on the other hand, to have $ \mathrm{Re}(X) \le -2 $, $ Y+2Z\le 0 $ is needed. By using parameters in~\eqref{check}, we get $ Y+2Z = 2d_0^2 + \mathrm{\bar W_2} \omega^2 +  \mathrm{\bar W_4} \omega^4 +  \mathrm{\bar W_6} \omega^6 +  \mathrm{\bar W_8} \omega^8 +  \mathrm{\bar W_{10}} \omega^{10} \le 0$, where $ \mathrm{\bar W_2} = d_0(2d_2+n_2) +2d_2d_0 + d_1(2d_1+n_1), 
													\mathrm{\bar W_4} = d_0(2d_4+n_4)+ d_1(2d_3+n_3)  + d_2(2d_2+n_2)+2d_0d_4 + d_3(2d_1+n_1),
													\mathrm{\bar W_6} = d_0n_6+ d_1n_5  + d_2(2d_4+n_4) + d_3(2d_3+n_3) + d_4(2d_2+n_2),
													\mathrm{\bar W_8} = d_2n_6+ d_3n_5  + d_4(2d_4+n_4),
													\mathrm{\bar W_{10}} = d_4n_6$. One can see as $ d_0 = k_1^2 >0 $ from~\eqref{k1_k2_k3}, even though we ask $ \{\mathrm{\bar W_2}, \mathrm{\bar W_4}, \mathrm{\bar W_6} \mathrm{\bar W_8}, \mathrm{\bar W_{10}}\} \le 0 $, it is still difficult to guarantee $ Y+2Z \le 0 $ for all $ \omega \in [0, \infty) $. This is the reason we  choose $ \mathrm{Re}(X) \ge 0 $ to have $ \left\|H(j\omega)\right\|_{\infty} \le 1 $.}
											\end{remark}
											
											\begin{remark}
												{Note} that {from Remark~\ref{reX_2},} the parameter setting condition~\eqref{parameter_setting_condition} is sufficient but not necessary for string stability condition~\eqref{string_stability_condition}. The advantage, however, is that designing controller parameters becomes much simpler  in both calculation and analysis.
											\end{remark}
											
											We recall that parameters  $ k_1,k_2,k_3 $ and $ \alpha \, (\bar \alpha = \alpha / \tau) $ remain to be determined. 
											In order to simplify the parameter setting, we design 
											\begin{align}\label{k1_k2_k3}
												k_1 = b^3\tau,\quad k_2 = 3b^2\tau,\quad k_3 = 3b\tau -1,\quad b>0 
											\end{align} 
											such that $ \lambda(A-BK)=-b, -b, -b $. 
											{This matrix $ K $  designing idea~\eqref{k1_k2_k3} delivers two advantages:}
											\begin{itemize}
												\item in this way, we  only need to design two scalar parameters only; namely, $ \alpha $ and $ b $;
												\item {designing $ b $ instead of $ K $ provides the parameter tuning convenience for dealing with the peaking effect of our proposed observer.}
											\end{itemize}

											\vspace{-0.3cm}
											
											\subsection{Heuristic searching algorithm for designing $ \alpha $ and $ b $}\label{sec_coupled_parameters_design}
											{We first present the following two observations for designing $ b $.}
											\begin{itemize}
												\item[I] {Polyak and Smirnov, in~\cite[Sec. 2]{polyak2016large} have shown that with all eigenvalues equal real $ -b<0 $, the large deviation effect is present both for $ b $ large and $ b $ small. However, the situation is different: for $ b $ large we have peaking effect in the initial period of time, while for $ b $ small the trajectory itself has large values and it happens for time $ t $ large enough.}
												{Authors in~\cite{smirnov2009advances} also proved that the peaking effect at the beginning of stabilization occurs with $ b \ll 1 $ and $ b\gg 1 $.
												}
												\item[II] Note that the larger the value of $ b $, the faster the convergence speed of the platoon, as $ -b $ is the eigenvalue of $ A-BK $ in the leader-following platooning  error dynamics~\eqref{x_tilde2}.
											\end{itemize}
											{To sum up, the first rule of designing $ b $ is that $ b $ should be large, but not too large, e.g., $ b\gg 1 $.}
											
											{Therefore,}
											from~\eqref{k1_k2_k3}, as the value of the engine time constant $ \tau $ is usually around $ 0.5 $ in real vehicles, based on $ b >0 $, we choose designing 
											$ b > 1/(3\tau) $. 
											In this way, we get $ \{k_1,k_2,k_3\} > 0 $.
											
											To give more details about the parameter setting condition~\eqref{parameter_setting_condition}, it is easy to calculate that 
											\vspace{-0.2cm}
											\begin{subequations}
												\begin{align}
													Y =& \mathrm{W_2} \omega^2 +  \mathrm{W_4} \omega^4 +  \mathrm{W_6} \omega^6 +  \mathrm{W_8} \omega^8 +  \mathrm{W_{10}} \omega^{10} \ge 0, \label{detail_parameter_setting_condition}\\
													\mathrm{W_2} =& d_0n_2 + d_1n_1, \label{W2}\\
													\mathrm{W_4} =& d_0n_4+ d_1n_3  + d_2n_2 + d_3n_1,\label{W4}\\
													\mathrm{W_6} =& d_0n_6+ d_1n_5  + d_2n_4 + d_3n_3 + d_4n_2,\\
													\mathrm{W_8} =& d_2n_6+ d_3n_5  + d_4n_4,\label{W8}\\
													\mathrm{W_{10}} =& d_4n_6. \label{W10}
												\end{align}
											\end{subequations}
											{From the definition of $ d_4 $ in~\eqref{d_4} and $ n_6 $ in~\eqref{n_6} with $ k_3 >0 $, we get $ \mathrm{W_{10}} <0 $.}
											In order to satisfy condition~\eqref{parameter_setting_condition}/\eqref{detail_parameter_setting_condition}, it is {sufficient} to have  
											\vspace{-0.2cm}
											\begin{subequations}\label{W_condition}
												\begin{align}
													\mathrm{W_2}, \mathrm{W_4}, \mathrm{W_6} &\ge 0,\label{W_condition1} \\
													\mathrm{W_6}+\mathrm{W_8} \omega^2 + \mathrm{W_{10}} (\omega^2)^2 &\ge 0 ,~\omega \in [0, \omega_0], \label{W_condition2}
												\end{align}
											\end{subequations}
											where $\omega_0$ is an upper bound on $\omega$ for which inequality \eqref{W_condition2} holds. 
											Note that \eqref{W_condition} is {a} sufficient but not {a} necessary condition for~\eqref{parameter_setting_condition}/\eqref{detail_parameter_setting_condition}. 
											As a consequence, {for} deciding the signs of $ \mathrm{W_2}, \mathrm{W_4}, \mathrm{W_6} $, {the variables $\alpha$ and $b$ should be chosen appropriately}.
											Based on {the fact that} $ \{k_1,k_2,k_3\} > 0 $, and since  $ r\ge 1 $, $ h>0$, and $\alpha >0 $, {then} from~\eqref{parameters}, it is obvious that $ \{n_1,n_4,n_5,d_0, d_4\} >0 $ and $ \{n_2,n_3,n_6\} <0 $. Therefore, the signs of $ d_1, d_2 $ and $ d_3 $ remain to be decided when designing $ \alpha $ and $ b $ {to satisfy the transformed string stability condition~\eqref{detail_parameter_setting_condition}}.

											We assume the predecessor number $ r $ and platoon time headway $ h $ are predefined {and summarize our $ \alpha $ and $ b $ designing mechanism as follows (the details of steps i and ii are in the Appendix~\ref{appendix_heuristic_searching}).}
											\begin{enumerate}
												\item [i.] 
												The main rule:
												\vspace{-0.2cm}
												\begin{subequations}
													\begin{align}
														&\frac{4\alpha(r-1)}{9\tau ^2} + \frac{8}{9\tau} \le b < \frac{6}{h}, \label{parameter_b}\\
														&\text{$ b $ should be large, but not too large.} \label{parameter_b1}
													\end{align}
												\end{subequations}
												\item [ii.] 
												The complementary rule:
												\vspace{-0.2cm}
												\begin{equation}\label{parameter_b_2}
													3b\tau^2(hb-5) + 2\tau - \alpha <0, \alpha >0.
												\end{equation}
												\item [iii.] 
												After $ \alpha $ and $ b $ are designed from the above two rules,
												by verifying $ \{ \mathrm{W_{4}}, \mathrm{W_{6}}, \mathrm{W_{8}} \} >0 $ and  $ \mathrm{W_8} + \mathrm{W_{10}} \omega^2 \ge 0 $ for \textit{some value range $ \omega $\footnote{It is easy to deduce that it is not possible to verify $ \mathrm{W_8} + \mathrm{W_{10}} \omega^2 \ge 0 $ for all values of $ \omega $.}} (here, we recommend $ \omega \in [0, 100 \text{rad/s}]$, 
												{i.e., $\omega_0=100 \text{rad/s}$}~\footnote{The dominating frequency range on the body ({spring} mass) is for passengers cars approximately 1-2 Hz; considering different road and speed conditions, the peak frequency is around 10 Hz (1Hz = 2$ \pi $ rad/s)~\cite{barbosa2012vehicle}. Here, we choose [0, 16 Hz].}), it  still does not guarantee condition~\eqref{detail_parameter_setting_condition} completely for \textit{all values of $ \omega  {\in [0,\omega_0]} $}. 
												\item [iv.] {A} Bode plot of the original string stability condition~\eqref{string_stability_condition} will be made to verify the designed parameters $ \alpha, b $.  If  Bode plot shows $ \left\|H(j \omega)\right\|_{\infty} \le 1 $ for all values of $ \omega {\in [0,\omega_1], \omega_1 \gg \omega_0} $, this mechanism is completed. Otherwise, we redesign $ \alpha $ and $ b $ from the main and complementary rules and repeat steps $ \text{iii, iv} $.
												\item [v.] If we cannot exit step $ \text{iv} $, it means for current values of $ r $ and $ h $, there is no solution from this mechanism. Change the values of $ r, h $ based on the platooning requirement and repeat steps $ \text{iii, iv} $. If step $ \text{iv} $ still cannot be exited, { then the proposed mechanism does not give feasible solutions} for the current platooning requirement, i.e.,  maximum predecessor number $ \bar r $, largest time headway $ \bar h $.
											\end{enumerate}
											\begin{remark}\label{alg_remark}
												Given a fixed time headway $ h $, the main rule~\eqref{parameter_b} and the complimentary rule~\eqref{parameter_b_2} for setting parameters $ \alpha $ and $ b $ are neither sufficient nor necessary to guarantee the original string stability condition~\eqref{string_stability_condition}. However, we can regard the above rules as a \textit{heuristic searching algorithm} for finding appropriate parameters $ \alpha $ and $ b $.
												Note that rules~\eqref{parameter_b} and~\eqref{parameter_b_2} are only necessary, but not sufficient conditions for $ \mathrm{W_{2}} \ge 0 $ and $ \mathrm{W_{4}} \ge 0 $, respectively. Therefore, there exist scenarios that with the designed $ \alpha $ and $ b $ from rules~\eqref{parameter_b} and~\eqref{parameter_b_2}, we do not have \eqref{W_condition1}; but we can use the Bode plot of the original string stability condition~\eqref{string_stability_condition} in step iv as a final measurement to verify $ \alpha $ and $ b $. If we have $ \left\|H(j \omega)\right\|_{\infty} \le 1 $ for all values of $ \omega \in [0,\omega_1] $, it means the designed $ \alpha $ and $ b $ are still available. 
											\end{remark}
											
											In practice, the platooning time headway $ h $ has a maximum value {$ \bar h $} from the {vehicle platooning} requirement {(too large value of $ h $ is not meaningful for the platoon in increasing road throughput and reducing travel time)}. In the following subsection, an algorithm is proposed to obtain the minimum acceptable value of $ h $.
											\vspace{-0.3cm}
											\subsection{Parameter $ \alpha $ design simplification and time headway minimization}\label{sec_h_minimization_design}
											\begin{algorithm*}[t]
												\caption{Time headway $ h $ minimization and controller parameter $  b $ design}\label{alg} 
												\begin{algorithmic}[1]
													\STATE \textbf{Input:} engine time constant $ \tau $, {largest predecessor number $  \bar r $, largest time headway $\bar h $ from the platooning specifications}, {frequency $\omega_0, \omega_1 \gg \omega_0$, parameter tuning maximum iteration number $ k_{\max} $},   tolerance $ \text{TOL} = 0.001 $. 
													\STATE {\textbf{Initialization:}} $k_1=0$, $ r = \bar r$, $h = \bar h$, $ h_{\mathrm{up}} = \bar h$,  $h_{\mathrm{lo}} = 0$, $h_{\mathrm{previous}} =0$, $\alpha = 2\tau$, $ \underline{b} = \frac{4\alpha(r-1)}{9\tau ^2} + \frac{8}{9\tau} $, $\mathrm{flag}_1=0$, $\mathrm{flag}_2=0$.		 
													\WHILE { $\mathrm{flag}_1=0$ } 
													\STATE $ \bar b(h) = \frac{5}{h}, \delta (h) = \frac{\bar b(h) -\underline{b} }{ k_{\max}} $. \COMMENT{given $ h $, $ \bar b(h) $: the upper bound of $ b $; $ \delta (h) $: the increment. }
													\WHILE {$ k_1 \le k_{\max} $ AND $\mathrm{flag}_2=0$}
													\STATE  $ b \leftarrow \underline{b}+ k_1 \delta (h) $.
													\COMMENT{guarantee that $b$ satisfies~\eqref{parameter_b2}}
													\IF { $ \left\|H(j \omega)\right\|_{\infty} \le 1 $ from condition~\eqref{string_stability_condition} for  $ \omega \in [0, \omega_1] $ }  
													\STATE  $\mathrm{flag}_2 \leftarrow 1$.
													\COMMENT{$\mathrm{flag}_2=1$ means the current value of $ h $ is verified and available.}
													\ELSE {} 
													\STATE $k_1\leftarrow k_1+1$.
													\COMMENT{increase the value of $ k_1 $ to increase the value of $ b $ in step 6. }
													\ENDIF
													\ENDWHILE
													\IF { $\mathrm{flag}_2=1$	}
													\STATE  $ h_{\mathrm{up}} \leftarrow h $,  $h_{\mathrm{previous}} \leftarrow h$, $  h \leftarrow (h_{\mathrm{lo}} + h_{\mathrm{up}})/2 $.
													\COMMENT{bisection-like algorithm for designing a smaller value of $ h $.} 
													\IF  {$ |h-h_{\text{previous}}| \le \text{TOL} $}
													\STATE $\mathrm{flag}_1 \leftarrow 1$. \textbf{Output} $ h_{\text{previous}}, b $.  \COMMENT{$\mathrm{flag}_1=1$ means the minimum $ h $ is found. exit Algorithm~\ref{alg}. }
													\ELSE {} 
													\STATE $k_1\leftarrow 0$,  $\mathrm{flag}_2 \leftarrow 0$. 
													\COMMENT{the new designed $ h $ from step 14 is not minimum.}
													\ENDIF
													\ELSE {} 
													\STATE    $ h_{\mathrm{lo}} \leftarrow h $, $  h \leftarrow (h_{\mathrm{lo}} + h_{\mathrm{up}})/2 $, $ k_1\leftarrow 0 $.
													\COMMENT{bisection-like algorithm for designing a larger value of $ h $.}
													\ENDIF
													\ENDWHILE
												\end{algorithmic}
											\end{algorithm*}

											From the rules~\eqref{parameter_b} and \eqref{parameter_b_2}, one can see designing $ b $ is related to design{ing} $ \alpha $, i.e., they are coupled. From the format of~\eqref{parameter_b_2} and $ hb<6 $ in \eqref{parameter_b}, we can simply set $ \alpha = 2\tau $ and choose $ b $ to satisfy $ hb<5 $. In this way, based on  the main rule~\eqref{parameter_b}, we {simplify the conditions needed for} designing $ b $ as {follows:}
											\begin{equation}\label{parameter_b2}
												\frac{4\alpha(r-1)}{9\tau ^2} + \frac{8}{9\tau} \le b < \frac{5}{h}.
											\end{equation}
											{Now that $\alpha$ is set and $r$ and $\tau$ are constant, }one can see that designing $ b $ is only related to the time headway $ h $. 
											Thus, after choosing an appropriate $ h $, the value range for choosing $ b $ from \eqref{parameter_b2} is fixed. {More specifically, we set the upper bound $ \bar b(h)$ as $ \bar b(h) \coloneqq \frac{5}{h}$ and by fixing $\alpha$, we obtain a lower bound as well, which is independent of $h$, given by $ \underline{b} = \frac{4\alpha(r-1)}{9\tau ^2} + \frac{8}{9\tau} $. Additionally, we form the function $\delta (h) \coloneqq \frac{\bar b(h) -\underline{b} }{k_{\max}} $, which states the step size by which $b$ will be increased in the algorithm. The parameter $k_{\max}$ is the maximum number of steps for which the algorithm is changing $b$ until a feasible solution is found.}
											
											The above controller parameters designing mechanisms are taken into account and Algorithm~\ref{alg}, which is followed by every vehicle in the platoon, is proposed. 
											In fact, Algorithm~{\ref{alg}} consists of two coupled sub-algorithms. In the first one, we find controller parameter $b$ and in the second one, we find the minimum acceptable value of the time headway $h$. The value found for the time headway in the second sub-algorithm is verified in the first sub-algorithm, in order to ensure the string stability. Specifically, Algorithm~{\ref{alg}} makes use of the following ideas:
											\begin{list4}
												\item Use the  Bode plot of the original string stability condition~\eqref{string_stability_condition} to check whether the current value of $ h $ is available or not (while loop: steps 4-11).
												\item If yes (step 13), use the bisection-like algorithm to design a smaller value of $ h $ (step 14); if not, design a larger value (step 21) and rerun the while loop (steps 3-23).
												\item During steps 13-19, we check whether the available value of $ h $ coming from step 8 (now is $ h_{\text{previous}} $ in step 14) is minimum or not by step 15. If yes, we exit Algorithm~\ref{alg} and output $ h_{\text{previous}}, b $; if not, we rerun the while loop (steps 3-23) by the newly designed $ h $ in step 14.
											\end{list4}
											
											
											\begin{remark}
												The reason we do not use the  condition~\eqref{W_condition} and the  Bode plot of the original string stability condition~\eqref{string_stability_condition} to double check whether the current value of $ h $ is available in Algorithm~{\ref{alg}}, is that condition~\eqref{W_condition} may restrain the available value of $ b $ and that Bode plot is a stronger verifying measurement than condition~\eqref{W_condition} for parameter $ b $. 
											\end{remark}

											\section{{Controller convergence rate discussion}}\label{sec_convergence_analysis}
											{From the controller error dynamics~\eqref{x_tilde2} and  the observer estimating error dynamics \eqref{platooning_error_vehicle1}, \eqref{platooning_error_vehicle2} and \eqref{observer_estimating_error_dynamics}, one can see our controller convergence rate mainly depends on the eigenvalues of $(A-BK)$ with the eigenvalues of $(A-BK-{r_i} BL)$ have some influence as well. }
											
											{As $L = \alpha B^T $ in~\eqref{L}, $ K = \begin{bmatrix}
													k_1&k_2 &k_3 
												\end{bmatrix} $ in~\eqref{design_k} and $k_1 = b^3\tau, k_2 = 3b^2\tau, k_3 = 3b\tau -1$ in~\eqref{k1_k2_k3}, we have $A-BK =$
												$\begin{bmatrix}
													-b&0 &0 \\
													0&-b &0 \\
													0&0&-b
												\end{bmatrix}, A-BK-{r_i} BL = \begin{bmatrix}
													-b&0 &0 \\
													0&-b &0 \\
													0&0&-b-\frac{r_i \alpha}{\tau^2}
												\end{bmatrix}$.
												We see that parameter $b$ will dominate the convergence rate of our controller. More specifically, the larger the value of $b$ is, the faster the observer-based controller converges.}
											
											{To sum up, when designing our controller parameters $(\alpha, b)$, in addition to follow the algorithms in Sec. \ref{sec_coupled_parameters_design} and Sec. \ref{sec_h_minimization_design}, we also take the controller convergence rate into consideration.
											}
												%
												%
												%
												%
												\vspace{-0.2cm}
												\section{Examples}\label{sec:examples}
												\vspace{-0.1cm}
												
												\begin{table}[b]
													\centering
													\caption{Model and PID controller~\eqref{input_bian} parameters}
													\label{tab:Model Parameters}
													\begin{tabular}{|c|c|c|c|c|c|c|}
														\hline
														$N$   & $D$ &$\tau$  & $ h$ & $ r$ & $v_0$ &$a_0$  \\ \hline
														7 & 5 $m$ & 0.5 $ s $  & 0.198 $ s $ & 3& $20$ $m/s$ &$10$ $m/s^2$ \\ \hline		
													\end{tabular}
												\end{table}

												Here, 
												the internal stability and the string stability of our observer-based controller are validated. 
												Also,	the comparison of our controller (\eqref{hat_a_i} and \eqref{input}) in this paper and the following distributed PID controller~\eqref{input_bian} from~\cite{bian2019reducing} is demonstrated:
												\vspace{-0.2cm}
												\begin{align}
													u_i =& - \sum_{j=1}^{N} a_{ij}[k_p(\tilde{p}_i - \tilde{p}_j) + k_v(\tilde{v}_i - \tilde{v}_j) + k_a(\tilde{a}_i - \tilde{a}_j)] \nonumber\\& - a_{i0}(k_p\tilde{p}_i + k_v\tilde{v}_i  + k_a\tilde{a}_i ),\label{input_bian}
												\end{align}
												where $ k_p = 0.1, k_v = 1.67, k_a = 0.84 $.
												As such, we use the same platoon parameters and PID controller~\eqref{input_bian} parameters as in \cite[Fig. 3(c)]{bian2019reducing}, which
												is listed in Table \ref{tab:Model Parameters}. 
												{It is worth noting that this PID controller~\eqref{input_bian} is also widely used under MPF IFT in other platoon works; see, e.g., \cite{9801548,9301227}.}
												A platoon that consists of a leader and $7$ following vehicles is considered.
												$r=3$ means only vehicles $ 3, 2 $ and $ 1 $ can get information from the leader vehicle and vehicles $ 7, 6, 5, 4 $ are connected to three vehicles directly ahead, as shown in Fig.~\ref{fig_platoon}. 
												Vehicle $i$ starts at the point $-iD$
												and moves to follow the leader's {speed}. 
												All followers' speeds and accelerations are set as $ 0 $. {One can check that the initial condition of our platoon system is safe, i.e., the unsafe condition~\eqref{safety_condition} is not satisfied. }
												
\begin{figure*}[h]
													\centering
													~
													\begin{subfigure}[b]{0.3\textwidth}
														\includegraphics[width=\textwidth]{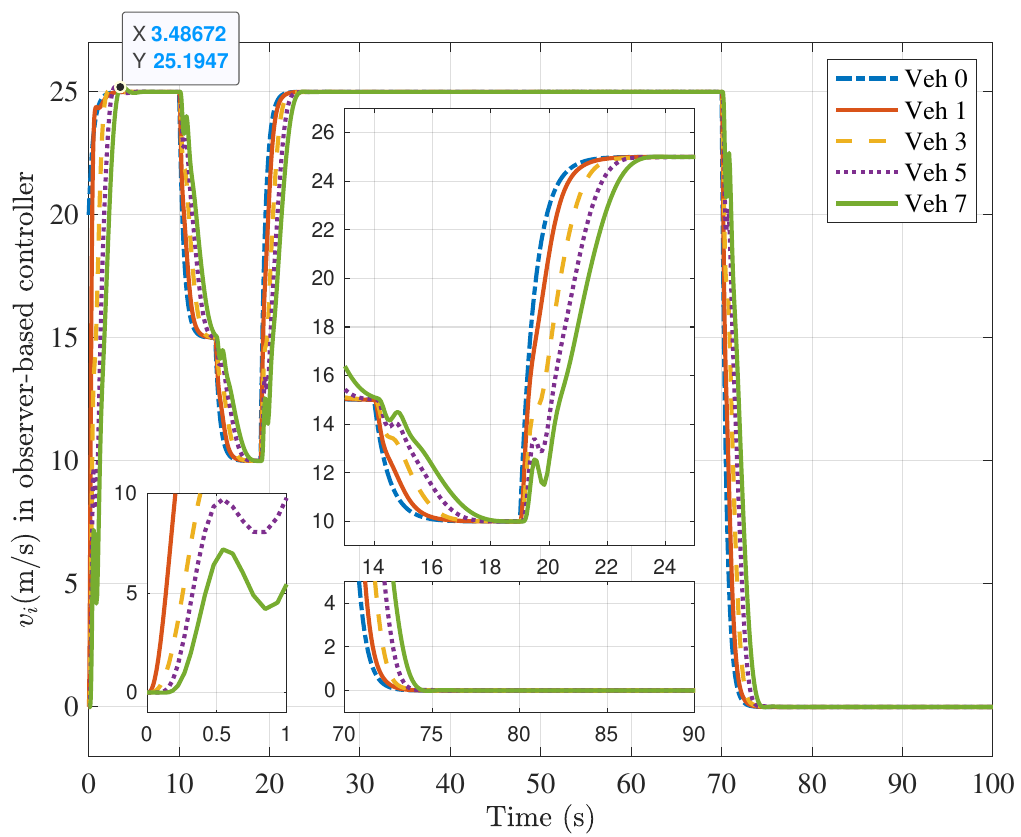}
														\caption{Observer-based controller (\eqref{hat_a_i} and \eqref{input}) with $ b=7.5 $.}
													\end{subfigure}
													~
													\begin{subfigure}[b]{0.3\textwidth}
														\includegraphics[width=\textwidth]{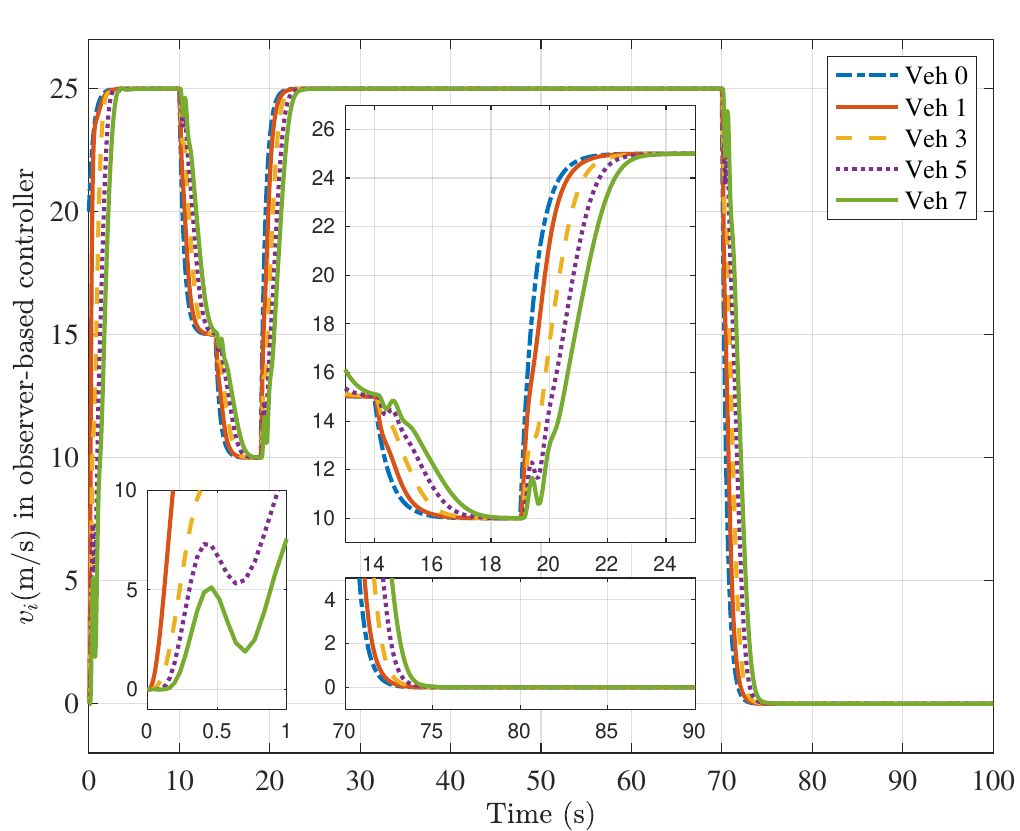} 
														\caption{ Observer-based controller (\eqref{hat_a_i} and \eqref{input}) with $ b=9 $. }
													\end{subfigure}
													~
													\begin{subfigure}[b]{0.3\textwidth}
														\includegraphics[width=\textwidth]{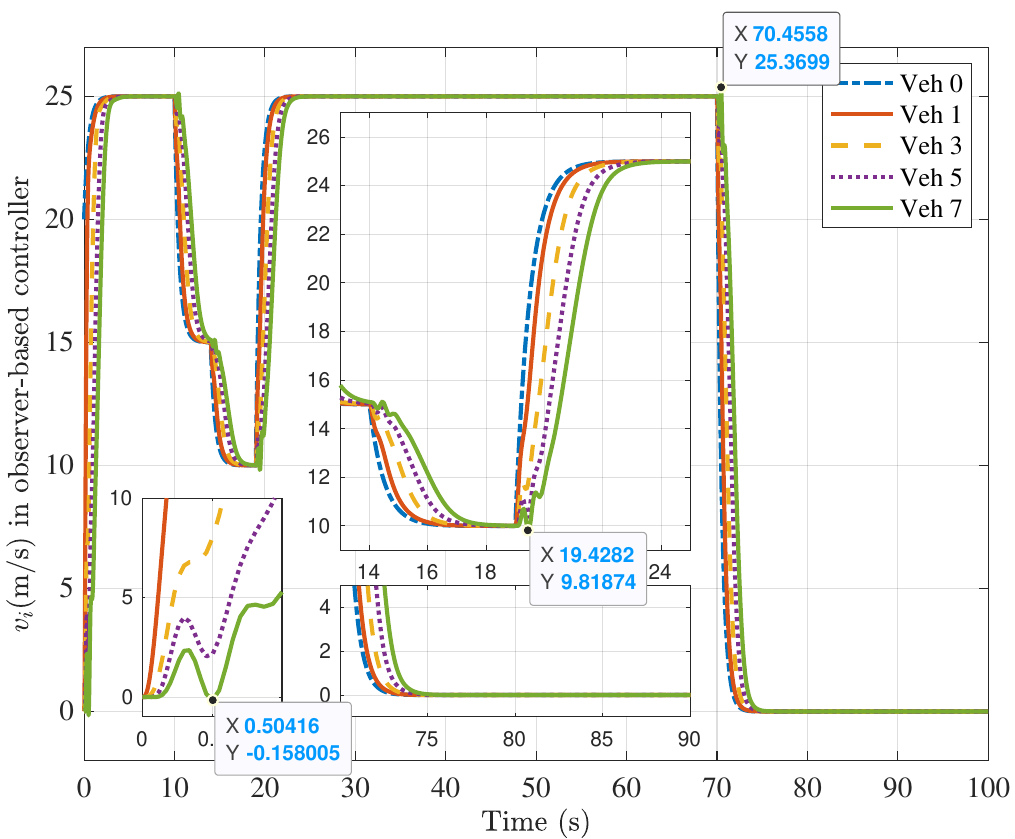} 
														\caption{Observer-based controller (\eqref{hat_a_i} and \eqref{input}) with $ b=12 $.}
													\end{subfigure}	
													\\
													\begin{subfigure}[b]{0.3\textwidth}
														\includegraphics[width=\textwidth]{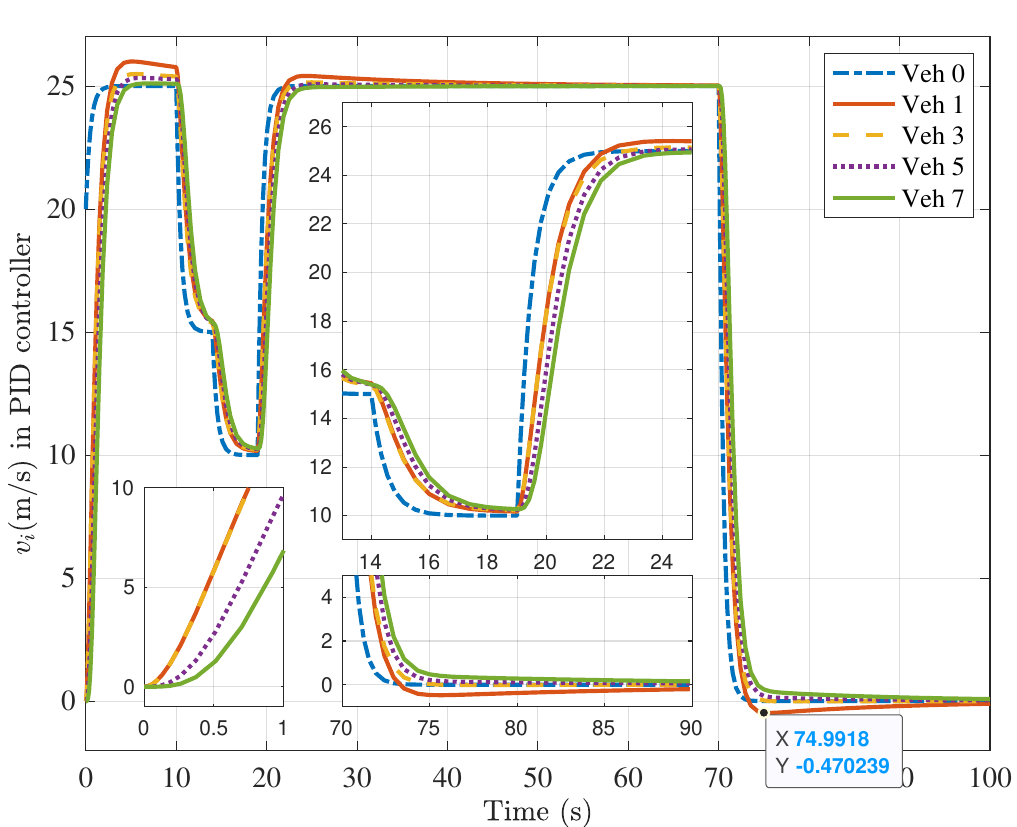} 
														\caption{Distributed PID controller~\eqref{input_bian} in~\cite{bian2019reducing}.}
													\end{subfigure}	
													~
													\begin{subfigure}[b]{0.3\textwidth}
														\includegraphics[width=\textwidth]{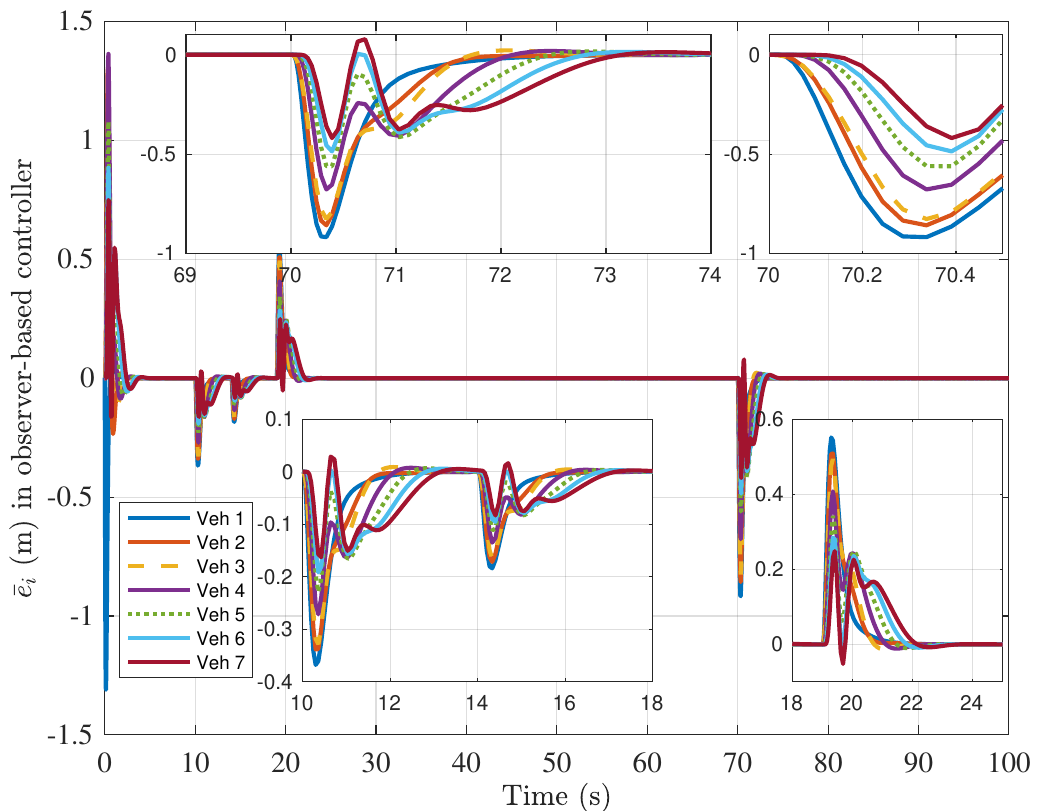}
														\caption{Observer-based controller  with $ b=9 $.}
													\end{subfigure}
													~
													\begin{subfigure}[b]{0.3\textwidth}
														\includegraphics[width=\textwidth]{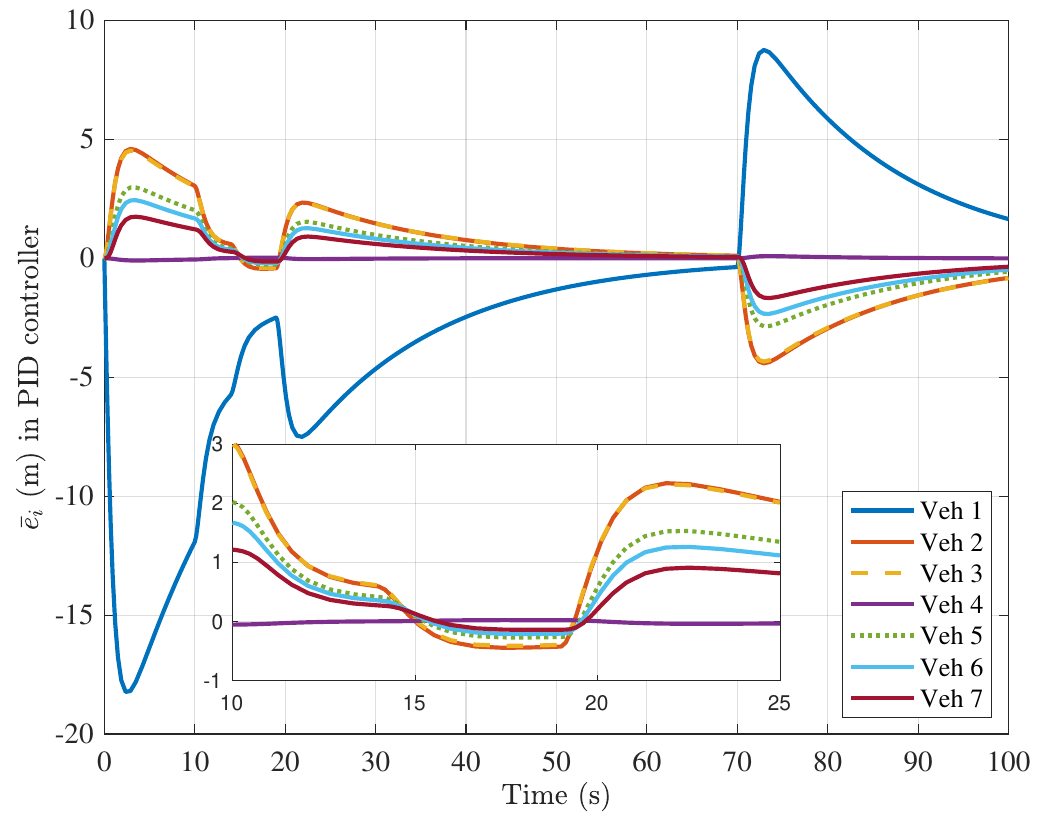} 
														\caption{Distributed PID controller~\eqref{input_bian} in~\cite{bian2019reducing}.}
													\end{subfigure}
													\caption{Performance comparison of followers tracking leader's speed in (i-iv) and
														comparison of platooning string stability performance with the predecessor-follower platooning spacing error $ \bar e_i $~\eqref{bar_e_i} in (v-vi).
													}
													\label{fig_compare}
												\end{figure*}

												To demonstrate the string stability, an external disturbance, $u_0(t) $ in the following is acted on the leader to have the leader speed profile as shown in Fig.~\ref{fig_compare} (iii) and (iv):
												\vspace{-0.2cm}
												\begin{equation}\label{leader_disturbance}
													u_0(t) =
													\begin{cases}
														0,  &  t\in [0s, 10s),  \\
														-20,  &  t\in [10s, 14s),  \\
														-30,  &  t\in [14s, 19s),  \\
														0,  &  t\in [19s, 70s),  \\
														{-50},  &  t\in[70s, 100s]. 
													\end{cases}
												\end{equation}
												{By designing the leader vehicle's input as~\eqref{leader_disturbance} in simulation, plus its dynamics~\eqref{leader_dynamics_pva} in theory, we generate the leader vehicle's speed profile as shown in Fig.~\ref{fig_compare} (iii) and (iv) similar to that in~\cite[Fig.~6]{lunze2018adaptive}.
												}  

												
												
												\vspace{-0.2cm}
												\subsection{Comparison with distributed PID controller~\eqref{input_bian} in~\cite{bian2019reducing} }\label{exampleA}

												%

												\begin{figure*}[t]
													\centering
													~
													\begin{subfigure}[b]{0.3\textwidth}
														\includegraphics[width=\textwidth]{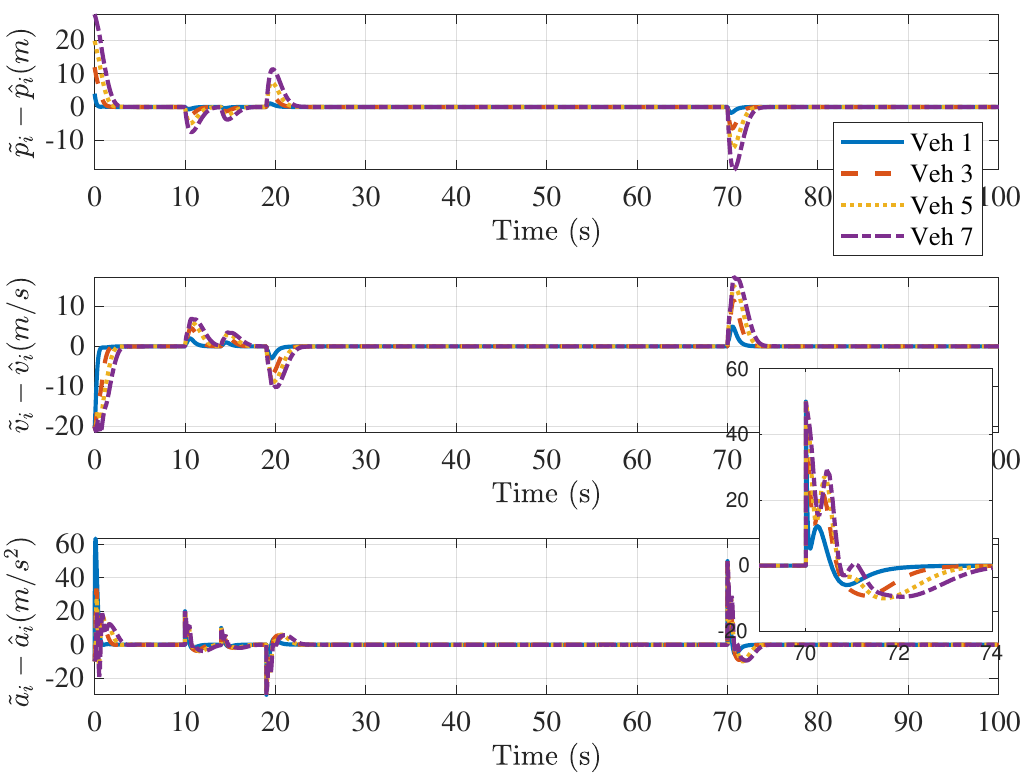}
														\caption{{Observer  error $ \xi_i $ in~\eqref{ESO_error_dynamics} between   output $ \hat{x}_{i} \coloneqq [\hat  p_i, \hat v_i, \hat a_i]^T $ and  real system unmeasurable states $ 
																\tilde x_i \coloneqq [\tilde p_i, \tilde v_i, \tilde a_i]^T 
																$.}
														}
													\end{subfigure}
													~
													\begin{subfigure}[b]{0.3\textwidth}
														\includegraphics[width=\textwidth]{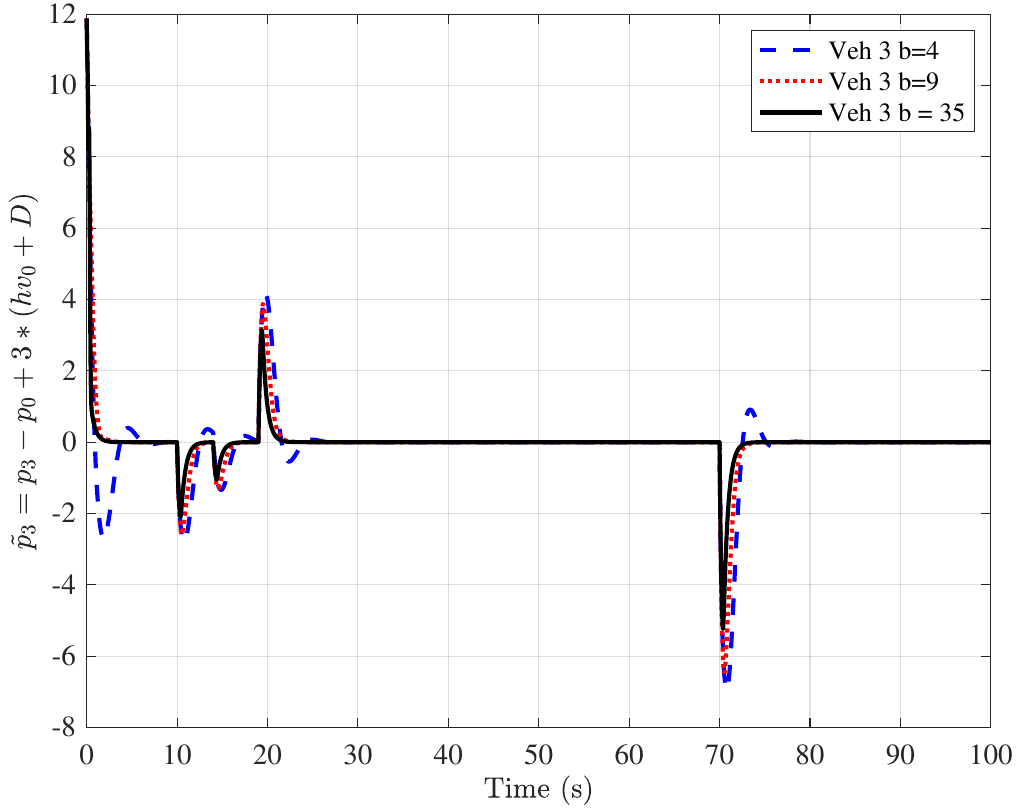} 
														\caption{{Controller convergence rate performance comparison for different values of $ b $ with $ \alpha = 1.5, h = 0.198s $.}}
													\end{subfigure}
													~
													\begin{subfigure}[b]{0.3\textwidth}
														\includegraphics[width=\textwidth]{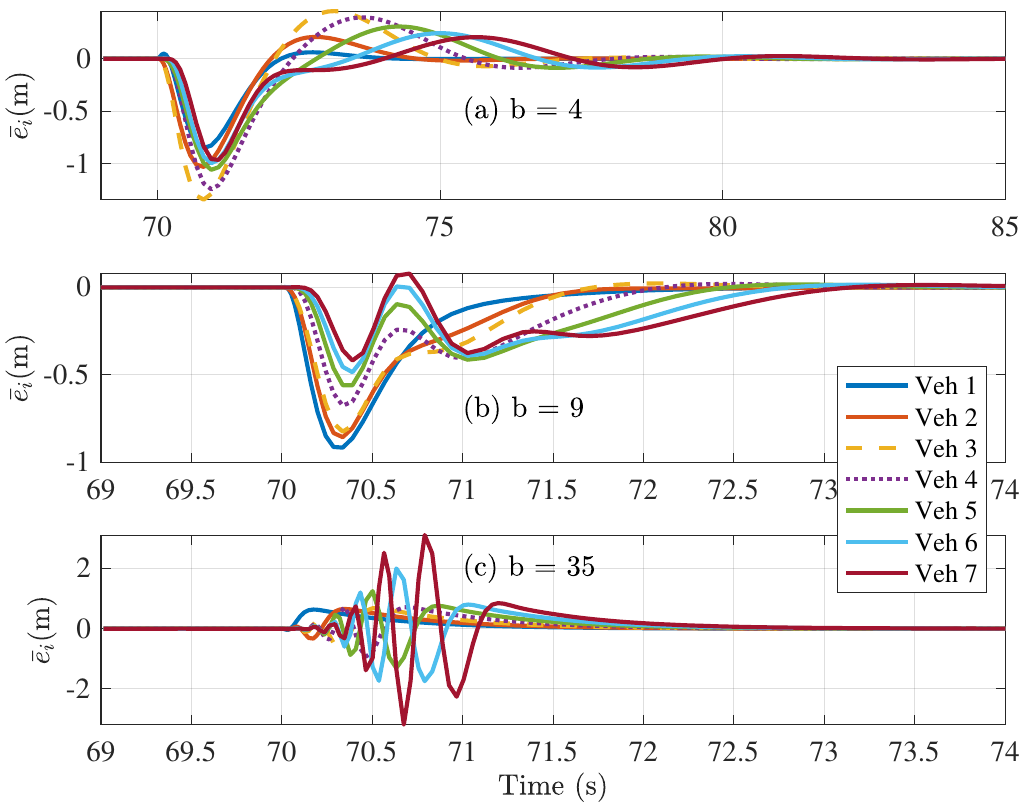} 
														\caption{Platooning string stability performance comparison for different values of $ b $ with $ \alpha = 1.5, h = 0.198s $.}
													\end{subfigure}
													\\
													\begin{subfigure}[b]{0.3\textwidth}
														\includegraphics[width=\textwidth]{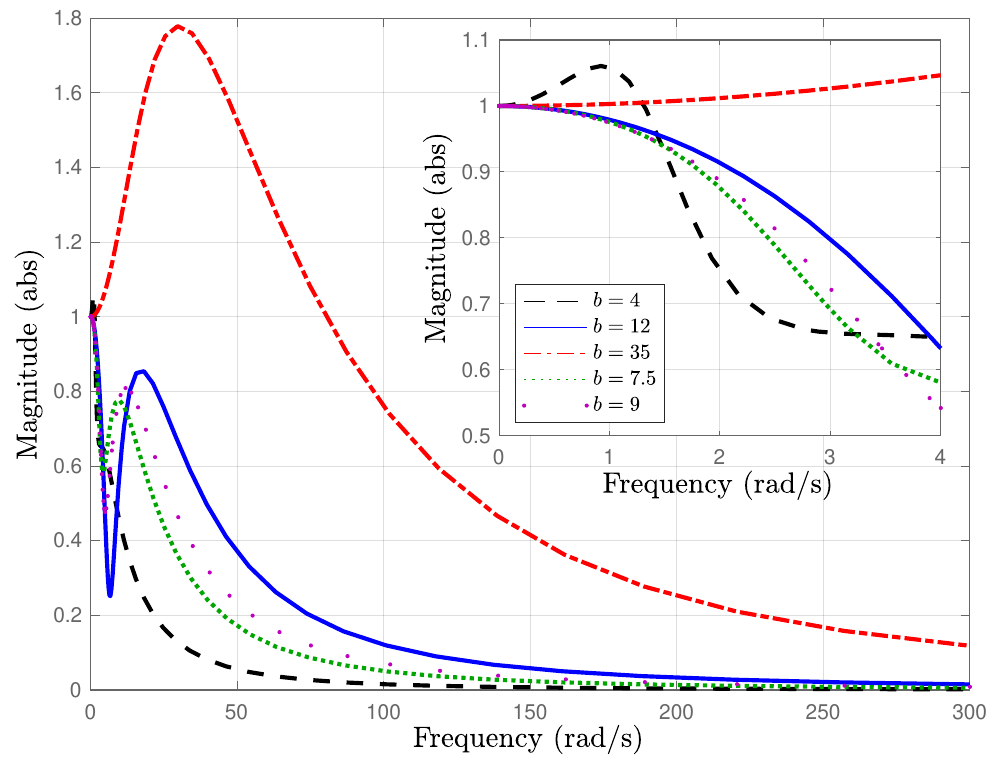}
														\caption{The Bode plot of the original string stability condition~\eqref{string_stability_condition} with different values of $ b $ with $ \alpha = 1.5, h = 0.198s $.}
													\end{subfigure}
													~ 
													\begin{subfigure}[b]{0.3\textwidth}
														\includegraphics[width=\textwidth]{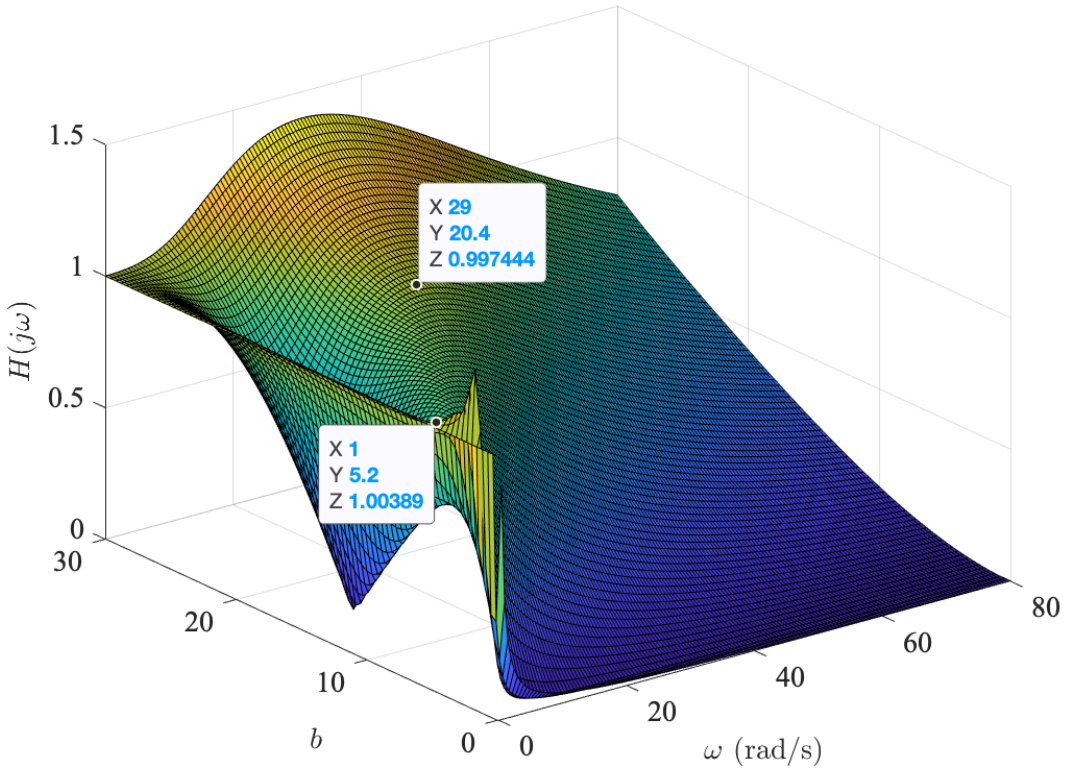} 
														\caption{The relationship between $ \left\|H(j \omega)\right\|_{\infty} $ and  $ b $ when $ \alpha = 1.5 $.}
													\end{subfigure}	
													~ 
													\begin{subfigure}[b]{0.3\textwidth}
														\includegraphics[width=\textwidth]{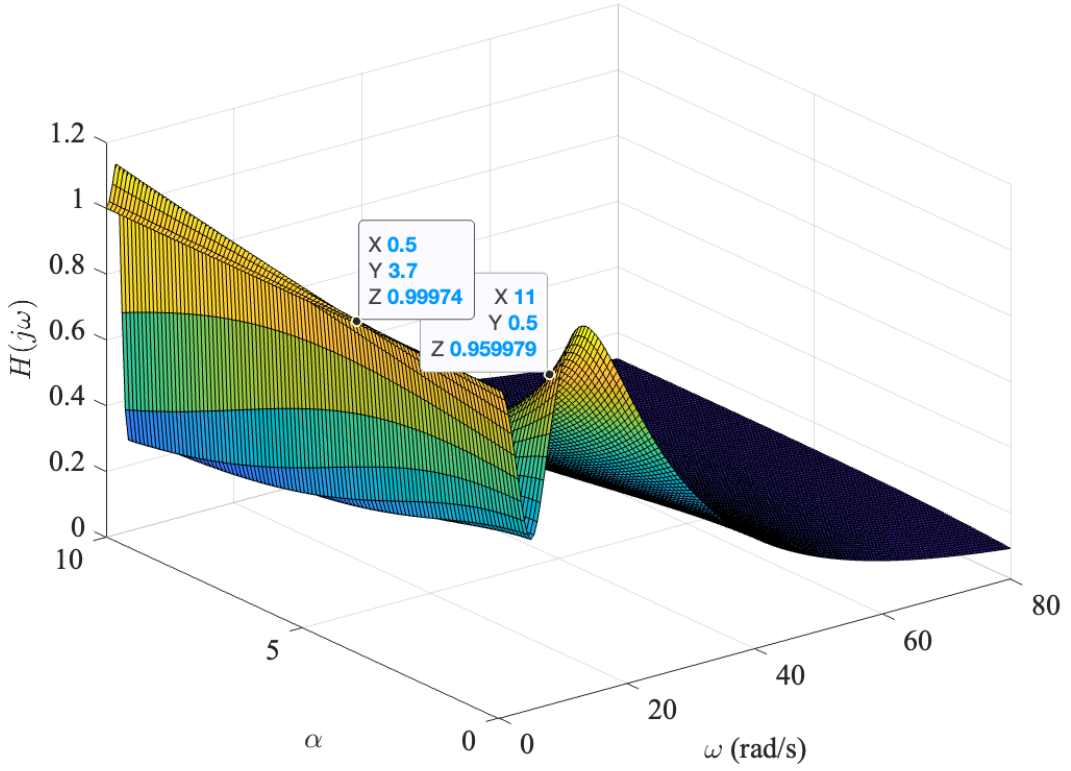} 
														\caption{The relationship between $ \left\|H(j \omega)\right\|_{\infty} $ and  $ \alpha $ when $ b = 9 $.}
													\end{subfigure}	
													\caption{Heuristic search algorithm validation for our controller parameters with fixed time headway $ h=0.198s $ in Sec.~\ref{exampleA1}.}
													\label{fig_platoonError_parameters}
												\end{figure*}

												Based on $ r=3 $ and $ \tau = 0.5 $,  from the complimentary rule~\eqref{parameter_b_2} in step ii, if we  design $ \alpha = 1.5 $, then, the main rule~\eqref{parameter_b} in step i  becomes $ 7.1 \le b < 30.3 $. Based on guidelines, {firstly,} we design
												$ b = 12 $. 
												We verify $ \{ \mathrm{W_{2}}, \mathrm{W_{4}}, \mathrm{W_{6}} \} \ge 0 $ in \eqref{W_condition1} as $ \mathrm{W_{2}} = 5.7773 e^{8}, \mathrm{W_{4}} = 1.0674 e^{8}, \mathrm{W_{6}}=7.3228 e^{6} $ and  $ \mathrm{W_6}+\mathrm{W_8} \omega^2 + \mathrm{W_{10}} (\omega^2)^2 \ge 0 $ in~\eqref{W_condition2} for $ \omega \in [0,\omega_0],\omega_0=3.9017  $ with $ \mathrm{W_{8}} = -1.7811 e^{4}, \mathrm{W_{10}} = -59.5 $. It is worth noting here that it does not mean the transformed string stability condition~\eqref{detail_parameter_setting_condition} is only available for $ \omega \in [0,\omega_0], $ as we also have $ \{\mathrm{W_{2}}, \mathrm{W_{4}} \} >0 $. However, the calculation for deciding the exact range for $ \omega $ with all $ \{\mathrm{W_{2}}, \mathrm{W_{4}} \} $ included would become very complicated as we can see in~\eqref{detail_parameter_setting_condition}. And this is the reason we use step iv to finally and formally confirm our designed feasible parameters $ \alpha = 1.5 $ and $ b = 12 $ as in the Fig.~\ref{fig_platoonError_parameters}  (iv) which means the string stability is guaranteed. 
												
												Fig.~\ref{fig_compare} (i-iv) demonstrate that:
												\begin{itemize}
													\item Follower vehicles using the observer-based controller {(Fig.~\ref{fig_compare} (i, ii and iii))} can track the leader's speed much faster and more smoothly.
													\item {When the leader's speed becomes very low, e.g., zero, the overshoot of the follower vehicle's speed may lead to negative speed, as shown for $ v_7 $ in Fig.~\ref{fig_compare} (iv), which may lead to a collision.}
												\end{itemize}
												
												\textbf{Peaking effect validation.}
												Note that for the presentation convenience, we only plot the speed profiles of vehicles $ 0, 1,3, 5, 7 $.
												{It is also interesting to find out that in Fig.~\ref{fig_compare} (iii) with $ b=12 $, during the initial time period, the speed of vehicle $ 7 $ is negative, which is due to the cause of the peaking effect in Sec.~\ref{sec_coupled_parameters_design}. Based on the main rule~\eqref{parameter_b1}, $ b=12 $ is large. Then, we design $ b=9 $ and we can see the peaking effect is avoided in Fig.~\ref{fig_compare} (ii). When we continue decrease the value of $ b $ until $ b=7.5 $ (still satisfies $ 7.1 \le b < 30.3 $), we see no peaking effect but overshoot in the beginning time period in Fig.~\ref{fig_compare} (i), which also validates our main rule~\eqref{parameter_b1}. Note that both $ b=9 $ and $ b=7.5 $ satisfy the string stability condition as shown in Fig.~\ref{fig_platoonError_parameters}  (iv).
												}

												Furthermore, Fig.~\ref{fig_compare} (v) and (vi) are used  to visualize the string stable performance.
												Specifically, Fig.~\ref{fig_compare} (v) shows clearly the internal and string stability of our proposed controller (with $ b=9 $). Compared with the PID controller in~\cite{bian2019reducing} from Fig.~\ref{fig_compare} (vi), one can see the error magnitude of our controller is much smaller, and its convergence speed is much quicker.

												{Fig.~\ref{fig_platoonError_parameters}  (i) describes the convergence to zero of the
													observer estimating  errors $ \xi_i $ in~\eqref{ESO_error_dynamics} between the observer ( with $ b=9 $) output $ \hat{x}_{i} \coloneqq [\hat  p_i, \hat v_i, \hat a_i]^T $ and  real system unmeasurable states: leader-following platooning tracking error $ 
													\tilde x_i \coloneqq [\tilde p_i, \tilde v_i, \tilde a_i]^T 
													$ in~\eqref{pva_error}. One can see  $ \xi_i $ converges to $ 0 $ one after one from vehicle $ 1 $ to vehicle $ 7 $, which validates Remark~\ref{remark_oneAfterOne}.}
												
												{\textbf{Controller convergence rate validation.} Fig.~\ref{fig_platoonError_parameters} (ii) takes vehicle 3 as an example to show the position error convergence rate comparison for different values of $b$ and demonstrates that with a larger value of $b$, our controller converges faster, which validates Sec.~\ref{sec_convergence_analysis}.}
												
												\vspace{-0.2cm}
												\subsection{Analysis of heuristic searching algorithm  in Sec.~\ref{sec_coupled_parameters_design}}\label{sec_herristic}
												
												{Since condition~\eqref{parameter_b} in the heuristic searching algorithm is only a necessary condition for $ \mathrm{W_2} \ge 0 $ in~\eqref{W_condition}, it is important to check the role that~\eqref{W_condition} plays in this  algorithm.}

												\subsubsection{Condition~\eqref{W_condition} is satisfied}\label{exampleA1}
												As we can see, the designed $ \alpha $ and $ b $ in Sec.~\ref{exampleA} shows~\eqref{W_condition} is satisfied.
												In order to further analyze and validate our proposed algorithm, we provide comparison simulations related to the different parameters of the proposed controller, e.g., by keeping $ \alpha = 1.5 $, we give two other choices of $ b $ ($ b=4 $  and $ b=35 $) which are out of the value range: $ 7.1 \le b < 30.3 $ from the main rule~\eqref{parameter_b}.
												The Bode plot of different values of parameter $ b $ in Fig.~\ref{fig_platoonError_parameters} (iv) demonstrates that the case of $ b=9 $ guarantees the string stability while the other two cases do not, which shows the effectiveness of our proposed heuristic searching algorithm (parameter setting mechanism) in~Sec.~\ref{sec_coupled_parameters_design}.

												More specifically, 
												one can see only Fig.~\ref{fig_platoonError_parameters} (iii) (b) demonstrates that the platoon is string stable as the error does not amplify from vehicle $ 1 $ to $ 7 $, while the scenarios of $ b=4 $ in (a) and $ b=35 $ in (c) are not, which is in accordance with the Bode plot in Fig.~\ref{fig_platoonError_parameters} (iv).

												\begin{figure*}[t]
													\centering
													~
													\begin{subfigure}[b]{0.3\textwidth}
														\includegraphics[width=\textwidth]{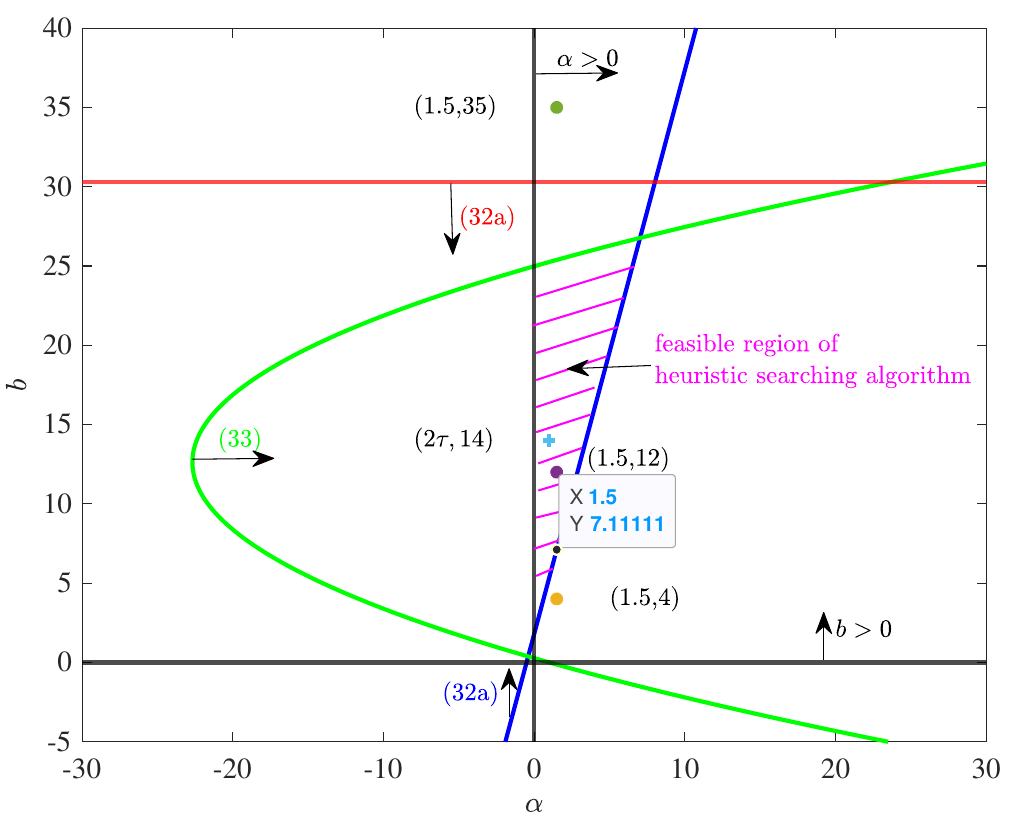}
														\caption{The feasible region for parameters $ \alpha $ and $ b $ given $ h=0.198 $. }
														\label{subfig:hr_greater_first_test_p}
													\end{subfigure}
													~
													\begin{subfigure}[b]{0.3\textwidth}
														\includegraphics[width=\textwidth]{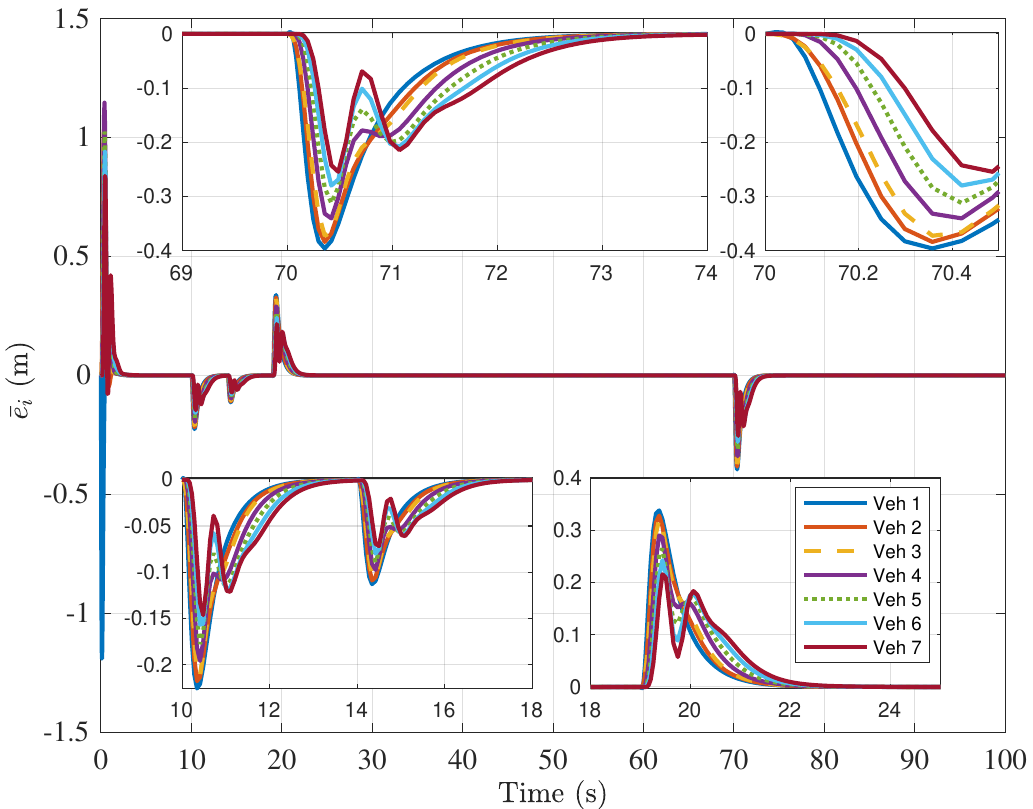} 
														\caption{Platooning string stability performance with $ \alpha = 2\tau, b=10, h=0.112 $.}
													\end{subfigure}	
													~ 
													\begin{subfigure}[b]{0.3\textwidth}
														\includegraphics[width=\textwidth]{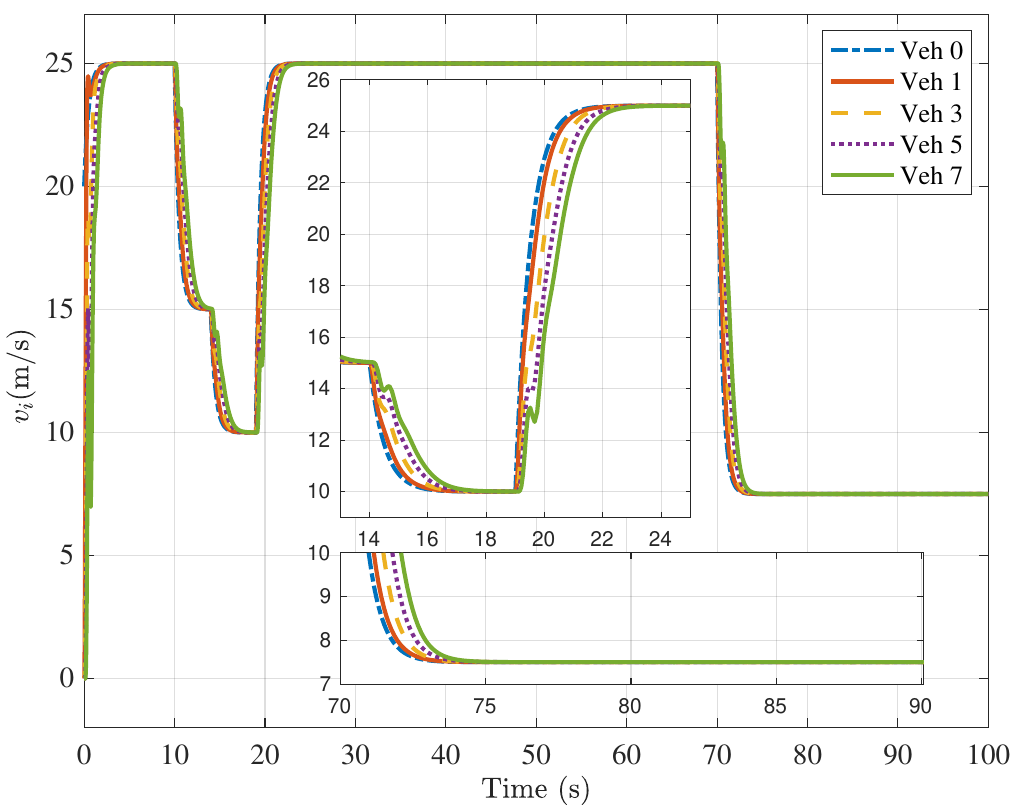} 
														\caption{Performance of followers tracking leader's speed with $ \alpha = 2\tau, b=10, h=0.112 $.}
													\end{subfigure}	
													\\
													\begin{subfigure}[b]{0.3\textwidth}
														\includegraphics[width=\textwidth]{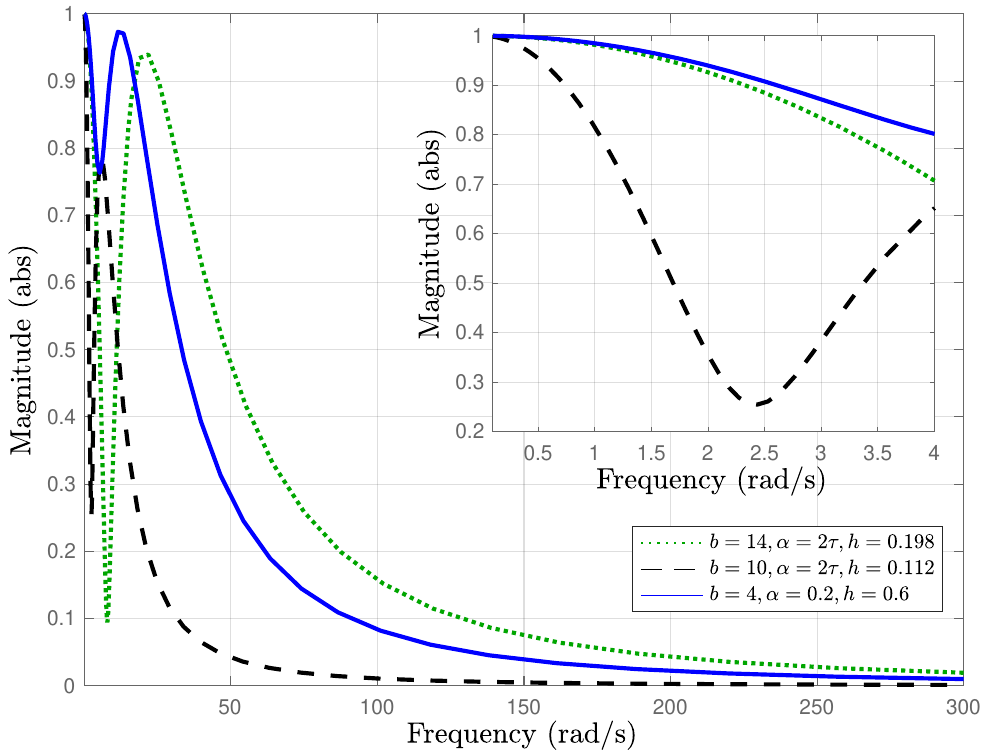} 
														\caption{The Bode plot  with different values of $ \alpha, b $ and $ h $.}
													\end{subfigure}
													~
													\begin{subfigure}[b]{0.3\textwidth}
														\includegraphics[width=\textwidth]{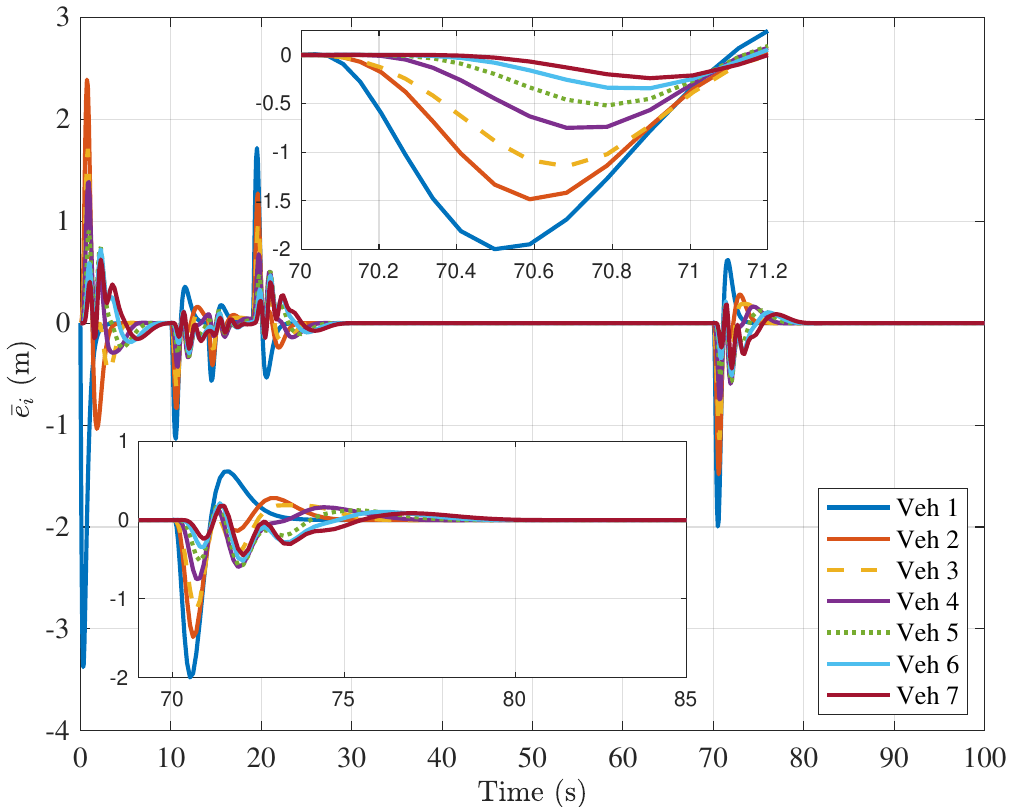} 
														\caption{Platooning string stability performance with $ h=0.6, \alpha = 0.2, b=4 $..}
													\end{subfigure}
													~ 
													\begin{subfigure}[b]{0.3\textwidth}
														\includegraphics[width=1\textwidth]{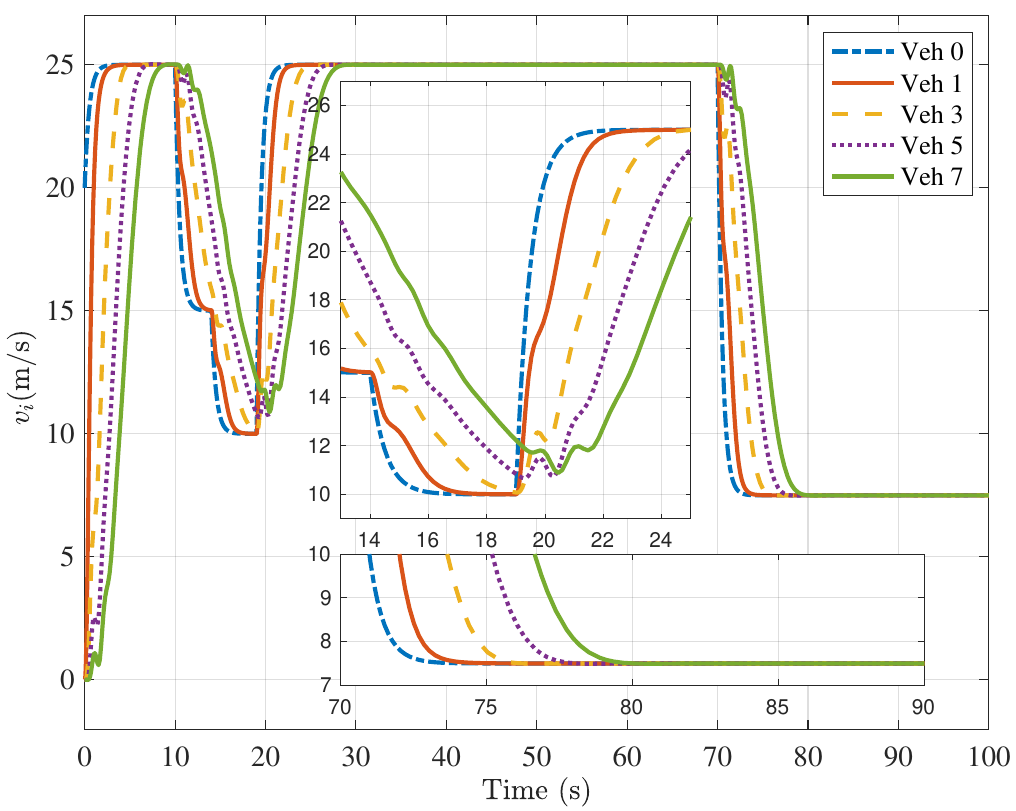} 
														\caption{Performance of followers tracking leader's speed with $ h=0.6, \alpha = 0.2, b=4 $..}
													\end{subfigure}	
													\caption{Feasible region for Sec.~\ref{exampleA1} and Algorithm~\ref{alg} for controller parameters with time headway minimization in Sec.~\ref{exampleB}.}
													\label{fig_h0112}
												\end{figure*}


												It is worth noting that all the errors converge to zero in  Fig.~\ref{fig_platoonError_parameters} (iii) (a-c) and that the error of vehicle $ 1 $ converges to zero first, then that of vehicle $ 2 $, and so on, and finally the error of the last vehicle (vehicle $ 7 $) converges to zero, which verifies Lemma~\ref{theorem_convergece_platooning} and Theorem~\ref{inter_string_theorem}.\\
												
												

												\textbf{Blind search method for controller parameters:}
												As we stated in Remark~\ref{alg_remark}, the parameter setting mechanism in~Sec.~\ref{sec_coupled_parameters_design} is neither sufficient nor necessary to guarantee the platoon string stability, and is rather a heuristic searching algorithm for possible solutions. 
												If we do not use the above mechanism, i.e., we do not have any information or guidelines about how to design $ \alpha $ and $ b $, then a blind search method {proposed in the conference version~\cite{jiang_platooning} of this paper} (i.e., guess the values of $ \alpha, b $ randomly and then tune them) can be adopted along with maybe many trials and errors.
												
												The blind search method is a numerical analysis method (e.g., \textit{meshgrid} and \textit{surf} functions in MATLAB) to get the relationship among $ \left\|H(j \omega)\right\|_{\infty} $ and parameters $ \alpha, b $ in \eqref{string_stability_condition}.
												To decide the relationship among $ \left\|H(j \omega)\right\|_{\infty} $ and parameters $ \alpha, b $, one solution is to use the command \textit{scatter3} in MATLAB to have a 4-D map which is not expressive to read. Instead, we choose to use commands \textit{meshgrid} and \textit{surf} to have 3-D maps which are more direct and clearer. 
												
												First, assume we set $ \alpha = 1.5 $ without many trials,  Fig.~\ref{fig_platoonError_parameters} (v) shows that $ 5.2 < b \le 20.4 $ is a decent choice for platooning string stability. One can see this range is different from the range  $ 7.1 \le b < 30.3 $ (or even $ 7.1 \le b < 25.25 $ considering the complimentary rule~\eqref{parameter_b_2} with $ b<5/h $), which verifies the above mechanism is neither sufficient nor necessary. However, they do have a quite large overlapped range $ 7.1 \le b \le 20.4 $ which accounts for $ 87.5\% $ of $ 5.2 < b \le 20.4 $ and $ 72.3\% $ of $ 7.1 \le b < 25.25 $, respectively, demonstrating the effectiveness of the above mechanism.

												Then, we choose $ b=9 $ and Fig.~\ref{fig_platoonError_parameters} (vi) shows $ 0.5 <\alpha \le 3.7 $ is acceptable.
												Therefore, the choice of $ (b=9, \alpha = 1.5) $ is good to satisfy the string stability condition~\eqref{string_stability_condition} given the time headway $ h $ is fixed.
												Fig.~\ref{fig_h0112} (i) shows the feasible region of our proposed controller parameter ($ \alpha $ and $ b $) design mechanism given $ h=0.198s $, which also validates the above analysis with $ (b=9, \alpha = 1.5) $ included.

												\subsubsection{Condition~\eqref{W_condition} is not satisfied}\label{exampleA2}
												
												Same as the above example, from the complimentary rule~\eqref{parameter_b_2} in step ii, if we design $ \alpha = 2\tau $, then, the main rule~\eqref{parameter_b} in step i  becomes $ 5.33 \le b < 30.3 $; we design
												$ b = 14 $. 
												Then, we have $ \mathrm{W_{2}} = -4.07 e^{9} <0, \mathrm{W_{4}} = 1.5136 e^{8} >0, \mathrm{W_{6}} = 1.046 e^{7} >0 $, which means condition \eqref{W_condition} is not satisfied.
												However, from Fig.~\ref{fig_h0112} (iv) we can see that by designing $  \alpha = 2\tau $ and $ b = 14 $, the string stability can be guaranteed.

												{This example verifies Remark~\ref{alg_remark} that  
													rules~\eqref{parameter_b} and~\eqref{parameter_b_2} are only necessary, but not sufficient for condition~\eqref{W_condition}. It confirms that the Bode plot of the original string stability condition~\eqref{string_stability_condition} is the formal and final measurement for string stability verification. 
													Also, this example demonstrates the effectiveness of our proposed rules~\eqref{parameter_b} and~\eqref{parameter_b_2} for guaranteeing the string stability of vehicle platooning.
												}
												
												\begin{remark}
													{The example in Sec.~\ref{sec_herristic} shows that given a fixed time headway $ h $, there exist multiple solutions of $ \alpha $ and $ b $ to guarantee the string stability of vehicle platooning. Sec.~\ref{exampleA} demonstrates the proposed observer-based controller controlling performance is better than the PID controller in~\cite{bian2019reducing}. Note that only two parameters, namely $ \alpha $ and $ b $, need to be tuned with the help of the feasible region, compared to at most three parameters in PID controller.}
												\end{remark}
												
												\subsection{Algorithm~\ref{alg} for controller parameters with time headway minimization in Sec.~\ref{sec_h_minimization_design}}\label{exampleB}
												
												Here, the same model parameters are used here, except that we set $ \alpha = 2\tau $ and use Algorithm~\ref{alg} to minimize $ h $ from $ h = \bar h =0.6s $. Finally, we arrive at $ h = 0.112s, b=10 $ from $ h = 0.6s, 0.3s, 0.15s, 0.075s, 0.112s $. Fig.~\ref{fig_h0112} (ii), (iii) and (iv) demonstrate the effectiveness and efficacy of proposed Algorithm~\ref{alg} for obtaining the minimum time headway $ h $.
												
												It is worth noting that the obtained minimum time headway $ h=0.112s $ is smaller than the obtained minimum time headway $ h=0.165s $ in \cite[Fig. 3(c)]{bian2019reducing} under the same condition.
												
												\begin{remark}
													One can observe that values of $ \alpha $ in the above two examples in Secs.~\ref{sec_herristic} and \ref{exampleB} are different. The reason is that Algorithm~\ref{alg} verified in  Sec.~\ref{exampleB} is by fixing $ \alpha = 2\tau  $ to minimize the value of $ h $. Actually, given a fixed $ h $, $ \alpha $ can have different values; and this is the heuristic searching algorithm in  Sec.~\ref{sec_herristic} demonstrates. To sum up, the objective of the heuristic searching algorithm is to design controller parameters $ \alpha $ and $ b $ given a fixed $ h $ while the objective of Algorithm~\ref{alg} is to minimize the value of $ h $ by fixing $ \alpha = 2\tau  $. To analyze the relationship among $ \alpha, b $ and $ h $ directly without fixing one value would be an interesting direction in the future.
												\end{remark}

												{\subsection{Controller performance  related to different values of $ h $}
													First, Fig.~\ref{fig_h0112} (iv) shows under the time headway $ h=0.6s $, the observer-based controller with $ \alpha = 0.2, b=4 $ is still string stable. 
													From Fig.~\ref{fig_h0112} (v), Fig.~\ref{fig_compare} (v) and Fig.~\ref{fig_h0112} (ii), one can see when the time headway becomes smaller and smaller, the string stability performance becomes smoother. It can be explained that the convergence rate of followers tracking leader's speed becomes faster without overshoot as shown in Fig.~\ref{fig_h0112} (vi), Fig.~\ref{fig_compare} (ii) and Fig.~\ref{fig_h0112} (iii). The reason is that with a larger time headway $ h $, when the leader's speed has a change (e.g., $ \Delta $), then the desired distance change (i.e.,  $ h \Delta $ ) between the follower and the leader is larger as we can see in~\eqref{deginition_ei}. Since there is no overshoot performance when the follower tracks the leader, then the follower takes longer time to track the leader's speed change. Note that even with a much larger time headway, from Fig.~\ref{fig_h0112} (vi) and Fig.~\ref{fig_compare} (iv), one can see the follower using observer-based controller tracks the leader faster than the one with PID controller.}
												
												\vspace{-0.1cm}
												%
												%
												%
												%
												\section{Conclusions and Future Directions}\label{sec:conclusions}
												
												\subsection{Conclusions}
												
												In this work, unlike the majority of literature on vehicle platooning, which assumes a constant speed leader, we study vehicle platooning control problem  with a leader whose {speed}  has a changing transient process with an exponentially converging behavior. Then, we design an observer-based controller under the directed {MPF} topology to improve the system’s performance when external disturbances occur.
												%
												%
												%
												The observer's matrix format is first proposed to guarantee the internal stability of the platoon system. Subsequently, by designing a specific observer parameter matrix, this observer turns out to have a third-order integrator dynamics (scalar format) which is utilized to derive the string stability conditions for designing observer and controller parameters. To deduce the string stability criterion, instead of calculating the derivatives of predecessor-follower spacing error directly which is also difficult, a new variable which is linked to that spacing error is proposed with its derivatives calculated instead until reaching the string stability criterion. 
												To design controller parameters from the above string stability criterion, a new parameter design mechanism given a fixed time headway is proposed to have a heuristic searching algorithm; furthermore, a bisection-like algorithm is incorporated into the {above} algorithm to  {obtain the minimum available}  value of the time headway by fixing one controller parameter.
												The validity and good platoon controlling performance of our proposed observer-based controller is demonstrated through comparison examples.

												\vspace{-0.2cm}
												\subsection{Future Directions}\label{future_work}
												
												{This work reveals a lot of opportunities to further enhance this idea towards its practical implementation. Some possible directions are discussed below.}
												\begin{list4}
													\item {As it is stated in~\cite{lunze2018adaptive},  string stability is only a necessary but not a sufficient condition for collision avoidance/safety; studying further the platoon's collision avoidance/safety conditions is of practical importance}.
													\item Investigating vehicle platooning with a leader whose speed is completely time-varying and has non-autonomous dynamics (i.e., non-zero input) is more realistic and challenging; part of ongoing work focuses on this problem.
													\item Another interesting direction is the study of the platoon system when the communication links are unreliable and they cause delays and packet losses.
												\end{list4}

												


												
												\appendices
												\section{Proof of Lemma~\ref{theorem_convergece_platooning}}
												\label{proof:lemma}
												First, we begin with the platooning convergence of vehicle $1$. Based on Eqs.~\eqref{ESO_error_dynamics}-\eqref{x_tilde2}, the observer error dynamics becomes
												\begin{align}\label{platooning_error_vehicle1}
													\dot{\xi}_{1}=(A-BK-BL)\xi_1 + B_1\Omega_1.
												\end{align}
												It is easy to have $ \lim_{t\rightarrow \infty} \xi_{1}(t) = 0 $ based on the condition that  $ (A-BK-BL) $ is Hurwitz ($ l_{ii}=1 $ as the leader is the only neighbor of vehicle $1$ in the platooning) and $ \lim_{t\rightarrow \infty} \Omega_1(t) = 0 $. Again from~\eqref{x_tilde2}, based on $ \lim_{t\rightarrow \infty} \xi_{1}(t) = 0 $,  $ \lim_{t\rightarrow \infty} \tilde{x}_{1}(t) = 0 $ can be proved based on $ A-BK $ is Hurwitz. Since $ \bar{x}_{1}=\tilde{x}_{1} + \begin{bmatrix}
													h\tilde v_{1}&
													0 & 
													0
												\end{bmatrix}^T $ {from~\eqref{relation_bar_tilde_x}}, we have $ \lim_{t\rightarrow \infty} \bar{x}_{1}(t) = 0 $.
												
												Then, for vehicle $ 2 $, similar to the previous calculation, we have 
												\begin{align}
													\dot{\xi}_{2}=&(A-BK)\xi_2 + B_1\Omega_2 - BL\sum_{j=1}^{2} l_{2j} \xi_j - \Pi_{1}\label{platooning_error_vehicle2} \\
													=&{(A-BK-r_{2}BL)\xi_2 + B_1\Omega_2 - BL l_{21} \xi_1 - \Pi_{1}. \nonumber}
												\end{align}
												{Note that $  r_i = l_{ii} , i \in \textbf{I}_1^N $ from Assumption~\ref{assump_mpf}.}
												From \eqref{modified_observer4}, one can see $ \Pi_{1} $ consists of $ \tilde p_1, \tilde v_1, \tilde a_1 $ which already all converge to zero because of $ \lim_{t\rightarrow \infty} \tilde{x}_{1}(t) = 0 $ for  vehicle $1$. In addition, based on  $ \lim_{t\rightarrow \infty} \xi_{1}(t) = 0 $  and $  \lim_{t\rightarrow \infty} \Omega_2(t) = 0 $, we get $ \lim_{t\rightarrow \infty} \xi_{2}(t) = 0 $. Furthermore, $ \lim_{t\rightarrow \infty} \tilde{x}_{2}(t) = 0 $ from~\eqref{x_tilde2}, {which also means $ \lim_{t\rightarrow \infty}\tilde v_{2}(t) = 0 $.
													In addition to $ \lim_{t\rightarrow \infty} \tilde{x}_{1}(t) = 0 $ for vehicle $ 1 $, we can have $ \lim_{t\rightarrow \infty} \bar{x}_{2}(t) = \lim_{t\rightarrow \infty} (\tilde{x}_{2}(t)-\tilde{x}_{1}(t)+ \begin{bmatrix}
														h\tilde v_{2}&
														0 & 
														0
													\end{bmatrix}^T) =0 $.}
												
												Finally, for vehicle $ i, i=3,4,\ldots, N $, we can get 
												{\small
													\begin{align}\label{observer_estimating_error_dynamics}
														\dot{\xi}_{i}
														=(A-BK-{r_i} BL)\xi_i + B_1\Omega_i  - BL \sum_{j=1}^{i-1} l_{ij} \xi_j  - \Pi_{i-1}. 
												\end{align}}
												From the MPF topology in Assumption~\ref{assump_mpf}, $ a_{ij}=1, l_{ij}=-1, j \in \textbf{I}_{i-r_i}^{i-1} 
												$ 
												and $ a_{ij}=0, l_{ij}=0, j \in \textbf{I}_{i-r_i-1}^{1}  
												$. In reality, $ \Pi_{i-1} $ in~\eqref{modified_observer4} consists of $  \tilde p_{i-1}, \tilde v_{i-1}, \tilde a_{i-1}  $ and $ \tilde v_{i-2}, \ldots, \tilde v_{i-r_i} $. From Eqs.~\eqref{platooning_error_vehicle1} and \eqref{platooning_error_vehicle2}, one can see for current vehicle $ i $, both $ \xi_j $ and $ \tilde{x}_j $ of its preceding vehicle $ j $ converge to zero asymptotically, i.e., $ \lim_{t\rightarrow \infty} \Pi_{i-1}(t) = 0  $. 
												{Thus, it is easy to have $ \lim_{t\rightarrow \infty} \xi_{i}(t) = 0 \Rightarrow\lim_{t\rightarrow \infty} \tilde{x}_{i}(t) = 0, i \in \textbf{I}_3^N $, which also includes $ \lim_{t\rightarrow \infty} \tilde{v}_{i}(t) = 0 $.
													Similarly, $ \lim_{t\rightarrow \infty} \bar{x}_{i}(t) = \lim_{t\rightarrow \infty} (\tilde{x}_{i}(t)-\tilde{x}_{i-1}(t)+ \begin{bmatrix}
														h\tilde v_{i}&
														0 & 
														0
													\end{bmatrix}^T) =0, i \in \textbf{I}_3^N $.} $\blacksquare $

												\section{Proof of Theorem~\ref{inter_string_theorem}}
												\label{proof:theorem:inter_string}
												For the internal stability, since the observer-based controller \eqref{hat_a_i} and \eqref{input} is transformed from controller \eqref{observer_hatx} and \eqref{input} by the Assumption~\ref{assump_mpf}, the internal stability proof is thus the same as the one in Lemma~\ref{theorem_convergece_platooning}.

												In order to prove string stability requirement: $\left\| \bar e_i \right\|_{\infty} \le \left\| \bar e_{i-1} \right\|_{\infty} $ where $ \bar e_i $ is defined in~\eqref{bar_e_i}, we propose another variable
												\begin{align}
													e_i =& p_i -p_{i-1} + h v_0+ D. \label{deginition_ei}
												\end{align}
												It is obvious from~\eqref{bar_e_i} that 
												\begin{equation}\label{relationship}
													\bar e_i = e_i+h {\tilde{v}_{i}}.
												\end{equation}
												
												In the following, we will construct the relationship between $ e_i $ and $ e_{i-1} $ instead. After that, the relationship between $ \bar e_i $ and $ \bar e_{i-1} $ can be built via~\eqref{relationship}.
												From vehicle dynamics~\eqref{platooning_dynamics} and input~\eqref{input}
												we  have 
												\begin{align}
													u_i = \tau \dot a_i + a_i =& \tau \dddot p_i + \ddot p_i = -k_1\hat p_i-k_2\hat v_i-k_3\hat a_i,\label{p_i}\\
													\tau \dddot p_{i-1} + \ddot p_{i-1}  =& -k_1\hat p_{i-1} -k_2\hat v_{i-1} -k_3\hat a_{i-1},\label{p_i1}\\
													\tau \dddot v_{0} + \ddot v_{0}  =& \tau \ddot a_{0} + \dot a_{0} = \tau \frac{-\dot a_{0}}{\tau} + \dot a_{0} =0.\label{dot_v_0}
												\end{align}
												Inspired from~\cite{bian2019reducing},
												by calculating \eqref{p_i} $ - $ \eqref{p_i1} $ + $ $ h\times $\eqref{dot_v_0}, we obtain
												\begin{align}
													\tau \dddot e_i + \ddot e_i 
													=& -k_1(\hat p_i - \hat p_{i-1})-k_2(\dot{\hat p} _i - \dot{\hat p}_{i-1})-k_3(\ddot{\hat p}_i - \ddot{\hat p}_{i-1}).\label{ei_pi}
												\end{align}
												Obviously, we need to calculate $ \hat p_i - \hat p_{i-1} $.
												From~\eqref{hat_a_i}, we get 
												\begin{align}
													\tau &(\dddot{\hat p}_i - \dddot{\hat p}_{i-1} ) = -k_1(\hat p_i - \hat p_{i-1})-k_2(\dot{\hat p} _i - \dot{\hat p}_{i-1})\nonumber\\&-(1+k_3)(\ddot{\hat p}_i - \ddot{\hat p}_{i-1})  +\sum_{l=1}^{r}\frac{\alpha}{\tau}(\ddot{\hat p}_{i-l} - \ddot{\hat p}_{i-l-1} ) \label{pi_minus1}\\&
													+ \sum_{l=1}^{r}\frac{\alpha}{\tau}[a_i-a_{i-1}-(a_{i-l}-a_{i-l-1})-(\ddot{\hat p}_i - \ddot{\hat p}_{i-1} )]\nonumber\\&
													+k_1[p_i -p_{i-1} - (p_{i-1} -p_{i-2}) + h(v_{i-1} -v_{i-2}) \nonumber\\&-(\hat p_i -\hat p_{i-1})]+k_2 [v_i- v_{i-1}-(v_{i-1} -v_{i-2})-(\dot{\hat p} _i - \dot{\hat p}_{i-1})]\nonumber\\&+k_3[a_i- a_{i-1}-(a_{i-1} -a_{i-2})-(\ddot{\hat p} _i - \ddot{\hat p}_{i-1})].\nonumber
												\end{align}
												From~\eqref{deginition_ei}, it is easy to get that
												$ p_i -p_{i-1} - (p_{i-1} -p_{i-2})= e_i - e_{i-1},
												v_i -v_{i-1}= \dot e_i - h\dot v_0=\dot e_i - h a_0, 
												a_i-a_{i-1}-(a_{i-l}-a_{i-l-1})=\ddot e_i -\ddot e_{i-l} $.
												Then, \eqref{pi_minus1} changes to 
												\begin{align}
													\tau &( \dddot{\hat p}_i - \dddot{\hat p}_{i-1} )\nonumber \\
													=& -2k_1(\hat p_i - \hat p_{i-1})-2k_2(\dot{\hat p} _i - \dot{\hat p}_{i-1})-(1+2k_3+ \frac{r\alpha}{\tau})\nonumber \\ & \times(\ddot{\hat p}_i - \ddot{\hat p}_{i-1})  +k_1e_i+k_2 \dot e_i +(k_3+ \frac{r\alpha}{\tau}) \ddot e_i  \label{pi_pi1} \\& -k_1e_{i-1} + (k_1h-k_2)\dot e_{i-1} -  k_3\ddot e_{i-1}  - \sum_{l=1}^{r}\frac{\alpha}{\tau}\ddot e_{i-l} + \Sigma_i, \nonumber
												\end{align}
												where $ \Sigma_i = \sum_{l=1}^{r}\frac{\alpha}{\tau}(\ddot{\hat p}_{i-l} - \ddot{\hat p}_{i-l-1} )-k_1 h^2  \ddot p_0 $.
												As we know that $ \lim_{t\rightarrow \infty} \hat{a}_i (t)= \lim_{t\rightarrow \infty}\tilde a_0(t)=0 $, so $ \lim_{t\rightarrow \infty} \Sigma_i(t)=0$.
												By setting all the initial conditions to be zero,  the Laplace transform of both sides of \eqref{ei_pi} and \eqref{pi_pi1} are
												\begin{align*}
													T_3E_i(s)=&- T_4(\hat P_i(s)-\hat P_{i-1}(s)), \\
													T_1(\hat P_i(s)-\hat P_{i-1}(s))=& T_2E_i(s)- \sum_{l=1}^{r}q_l(s) E_{i-l}(s) + \Sigma_i(s),
												\end{align*}
												where $ E_i(s), \hat P_i(s) $ are respectively the Laplace transformation of $ e_i(t), \hat p_i(t) $; {$ T_1, T_2, T_3, T_4, q_1(s) $ are defined in~\eqref{T1234} and $ q_l(s) = \bar \alpha s^2, l \in \textbf{I}_2^r $.}
												
												After simple mathematical manipulations, we come to
												\begin{align}
													E_i(s) =& \frac{q_1(s) T_4}{T_1T_3+ T_2T_4}E_{i-1}(s)\label{objective} \\& +\underbrace{s^2 \sum_{l=2}^{r} \frac{\frac{\alpha}{\tau^{2}} T_4}{T_1T_3+ T_2T_4}E_{i-l}(s) -\frac{ T_4\Sigma_i(s)}{T_1T_3+ T_2T_4} }_{\eqqcolon \varTheta_i(s)}. \nonumber
												\end{align}
												Recalling \eqref{relationship}, we have its Laplace transform as
												\begin{equation*}
													\bar E_i(s) = E_i(s)+h \tilde{V}_{i}(s) = E_i(s)+sh(P_{i}(s)-P_0(s)).
												\end{equation*}
												As a result, we finally have the string stability condition
												\begin{align}
													&\bar E_i(s) = \frac{q_1(s) T_4}{T_1T_3+ T_2T_4} \bar E_{i-1}(s)+ \varTheta_i(s) 
													\label{objective2} \\ &+
													\underbrace{sh( P_{i} (s)-P_0(s)) - \frac{shq_1(s) T_4}{T_1T_3+ T_2T_4}( P_{i-1} (s)-P_0(s))}_{\eqqcolon \Psi_i(s)}. \nonumber
												\end{align}
												The objective now is to design $ k_1, k_2, k_3, \alpha $ (or $ K, L $) such that $\left\| \bar E_i(j\omega) \right\|_{\infty} \le \left\| \bar E_{i-1}(j\omega) \right\|_{\infty} $.
												Note in~\eqref{objective} that $ T_4\Sigma_i(s) = T_4s^2 \Delta_i(s), \Delta_i(s) \eqqcolon  \sum_{l=1}^{r}\bar \alpha(\hat P_{i-l}(s) - \hat P_{i-l-1}(s ))+ k_1h^2 P_0(s)  $.
												Thus, we can see  $ \varTheta_i(s)  $ in~\eqref{objective} has two zeros located at the origin. Similarly, $ \Psi_i(s) $ has one zero located at the origin.
												Similar to \cite[Eq. (59)]{7134769}, if we have the following conditions at low frequencies 
												\begin{align}
													\left\| \frac{q_l(j\omega) T_4}{T_1T_3+ T_2T_4}E_{i-l}(j\omega) \right\|_{\infty} \ll 1, \omega &\rightarrow 0, l \in \textbf{I}_2^r, \label{extra terms2}\\
													\left\| -\frac{ T_4}{T_1T_3+ T_2T_4} \Sigma_i(j\omega) \right\|_{\infty} \ll 1, \omega &\rightarrow 0, \label{extra terms}\\
													\left\| \Psi_i(j\omega) \right\|_{\infty} \ll 1, \omega &\rightarrow 0, \label{psi}
												\end{align}
												then,  we can guarantee the string stability if the condition~\eqref{string_stability_condition} is satisfied~\cite{shaw2007string}.
												
												Now, we check whether  \eqref{extra terms2} and \eqref{extra terms} are satisfied.
												When $ \omega \rightarrow 0 $, we get $ T_1T_3 \rightarrow 0 \Rightarrow T_1T_3+ T_2T_4 = k_1^2 $; at the same time, as  $\omega \rightarrow 0 $, 
												$ q_l(j\omega) T_4 E_{i-l}(j\omega)~\approx~k_1 \bar \alpha  (-\omega^2)E_{i-l}(j\omega) $ and
												$ T_4\Sigma_i(j\omega)~\approx~k_1 (-\omega^2)\Delta_i(j\omega)$. Therefore, \eqref{extra terms2} and \eqref{extra terms} are satisfied.  
												
												Finally, for \eqref{psi}, as $ \Psi_i(j\omega) $ has only one zero located at the origin, it is a weaker inequality compared to \eqref{extra terms2} and \eqref{extra terms}. Nevertheless, at low frequencies \eqref{psi} is satisfied. $ \blacksquare $
												
												\begin{remark}\label{remark_x1}
													If we design $ \hat x_1 = \tilde x_1 $ for vehicle $ 1 $ instead of our proposed observer in~\eqref{observer_hatx1}, we can calculate the relationship between $ \bar e_2 $ and $ \bar e_1 $ directly. In this way, Eqs.~\eqref{p_i1} and \eqref{dot_v_0} will change respectively to $ 
													\tau \dddot p_{1} + \ddot p_{1}  = -k_1\tilde p_{1} -k_2\tilde v_{1} -k_3\tilde a_{1} $ and $ \tau \dddot v_{1} + \ddot v_{1}  = \tau \ddot a_{1} + \dot a_{1} = \dot u_1 =-k_1\dot{ \tilde p}_{1} -k_2 \dot{\tilde v} _{1} -k_3 \dot {\tilde a}_{1} $ with $ \bar e_1 = \tilde p_{1} $. Also, instead of calculating $ \dddot{\hat p}_i - \dddot{\hat p}_{i-1} $ in~\eqref{pi_minus1}, we calculate $ \dddot{\hat p}_2 $ only. Then, after some math manipulations, we conclude that for string stability, it is required that 
													\begin{equation}\label{H1_new_observer}
														\left\| H(j \omega) \right\|_{\infty} \coloneqq
														\left\| \frac{{\alpha \omega^2 T_4 + T_1 T_5} }{T_1T_3+ T_2T_4} \right\|_{\infty} \leq 1 ,
													\end{equation} where $ T_5 \coloneqq k_1 + (k_2-hk_1)s + (k_3-hk_2)s^2 -hk_3s^3 $. However, in this case, one can see when $ \omega = 0 $, $ \left\|H(j \omega)\right\|_{\infty} = 2 > 1 $. 
													Therefore, we design  observer $ \hat x_1  $ in~\eqref{observer_hatx1} for vehicle $ 1 $ instead of
													$ \hat x_1 = \tilde x_1 $.
												\end{remark}

												\section{The details of steps i and ii of proposed heuristic searching algorithm}\label{appendix_heuristic_searching}
												The following two steps is for  satisfying \eqref{W_condition} by dealing with $ d_1, d_2 $ and $ d_3 $ .
												\begin{itemize}
													\item [i.] \textit{$ d_1 >0, d_3 <0 $.} To have $ \mathrm{W_2} \ge 0 $ in~\eqref{W2}, {one should} have $ d_1 >0 $, i.e., $ h < 2k_2/k_1 $ or $ h< 6/b $. {In addition,} based on $ d_1 >0 $, one can confirm {that} $ d_3 <0 $ from~\eqref{d_3}. On the other hand, $ \mathrm{W_2} \ge 0 \Rightarrow f(h) \coloneqq k_1 h^2 -k_2 h + \bar \alpha (r-1) + 2 \le 0$.
													As $ h \in \mathbb{R} $,
													Based on the discriminant of the quadratic polynomial, we get $ (-k_2)^2 - 4k_1 (\bar \alpha (r-1) + 2) \ge 0 \Rightarrow b \ge \frac{4\alpha(r-1)}{9\tau ^2} + \frac{8}{9\tau} $. In addition to previous inequalities for parameter $ b $, we arrive at~\eqref{parameter_b}.
													\item [ii.] {\textit{$ d_2 <0 $ preferred.}} Since $ \mathrm{W_{10}} <0 $, even though if $ \{\mathrm{W_{2}}, \mathrm{W_{4}}, \mathrm{W_{6}} \} \ge 0 $, it is not obvious to guarantee the condition~\eqref{detail_parameter_setting_condition}.
													We switch to a different point of view. Going back to the original string stability condition~\eqref{string_stability_condition}, it is direct that $ \left\|H(j \omega)\right\|_{\infty} \rightarrow 1 $ when $ \omega \rightarrow 0 $ and $ \left\|H(j \omega)\right\|_{\infty} \rightarrow 0$ when $ \omega \rightarrow \infty $. It means if $ \left\|H(j \omega)\right\|_{\infty} >1 $, the value of $ \omega $ would not be very large. Therefore, when $ 0< \omega < 1 $, due to $ \mathrm{W_{2}}\ge 0 $, it has a very high possibility that $ \mathrm{W_2} \omega^2 + \mathrm{W_{8}} \omega^{8}+ \mathrm{W_{10}} \omega^{10} \ge 0 $ if $ b $ is properly selected from~\eqref{parameter_b}. When $ \omega \ge 1 $ and is not very large, by choosing $ b $ such that $ \{ \mathrm{W_{4}}, \mathrm{W_{6}} \} \ge 0 $,
													satisfying~\eqref{W_condition2} is usually available for a quite large value range  for  $ \omega $. Based on $ d_0n_4 >0, d_1n_3 <0, d_3n_1 <0 $ and $ n_2 <0 $, from the construction of $  \mathrm{W_{4}} $ in~\eqref{W4}, it is preferred to have $ d_2 <0 $, i.e., {from~\eqref{d_2}, we have~\eqref{parameter_b_2}.}
												\end{itemize}

												\section*{ACKNOWLEDGMENT}
												{The authors would like to thank anonymous reviewers and the Associate Editor for their important, enlightening and valuable comments.}

												

												\bibliographystyle{IEEEtran}
												\bibliography{mybibfile}

												%
												%
												%
												%
												\begin{IEEEbiography}[{\includegraphics[width=1in,height=1.25in,clip,keepaspectratio]{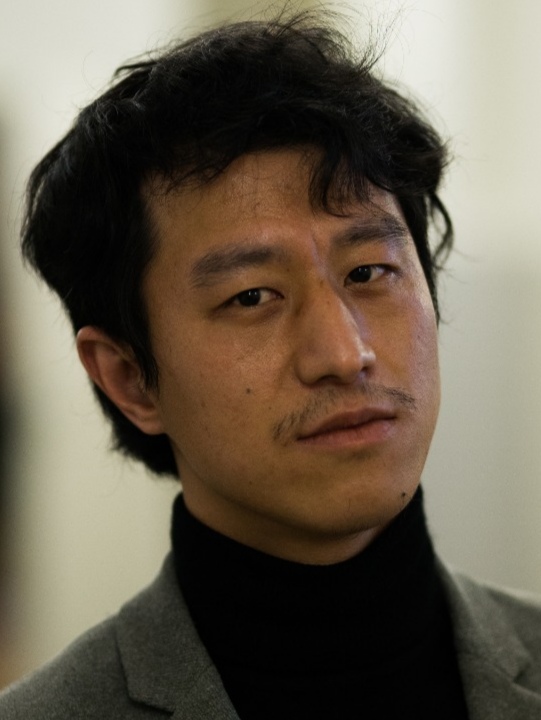}}]{Wei Jiang}
													received his B.S. degree in mechanical engineering and automation from Wuhan University of Technology, Wuhan, China, in 2011, and M.S. degree in automobile engineering from Beihang University, Beijing, China, in 2015 and Ph.D. degree at Automatic, Computer Engineering, Signal Processing and Images in CRIStAL, UMR CNRS 9189, Ecole Centrale de Lille, France, in 2018.
													He is a postdoctoral researcher at Aalto University, Finland since 2019. He was a Visiting Scholar at the University of Cyprus, Cyprus in October-November 2021 and KU Leuven, Belgium in April-May 2022. His research interests include vehicle platooning, distributed optimization/control,  learning algorithms and robotics.
												\end{IEEEbiography}
												
												\begin{IEEEbiography}
													[{\includegraphics[width=1in,height=1.25in,clip,keepaspectratio]{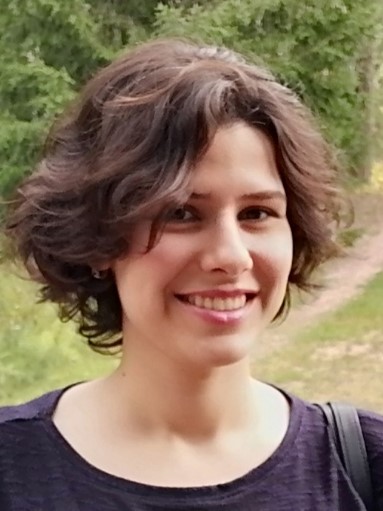}}]%
													{Elham Abolfazli}
													
													received the B.Sc. degree in electrical engineering from the Iran University of Science and Technology, Tehran, Iran, in 2012, and the M.Sc. degree in electrical and control engineering from the University of Tehran, Iran, in 2015. She is currently pursuing the Ph.D. degree with the Department of Electrical Engineering and Automation, School of Electrical Engineering, Aalto University, Finland. Her current research interests include networked control systems and connected vehicles.

												\end{IEEEbiography}
												
												\begin{IEEEbiography}[{\includegraphics[width=1in,height=1.25in,clip,keepaspectratio]{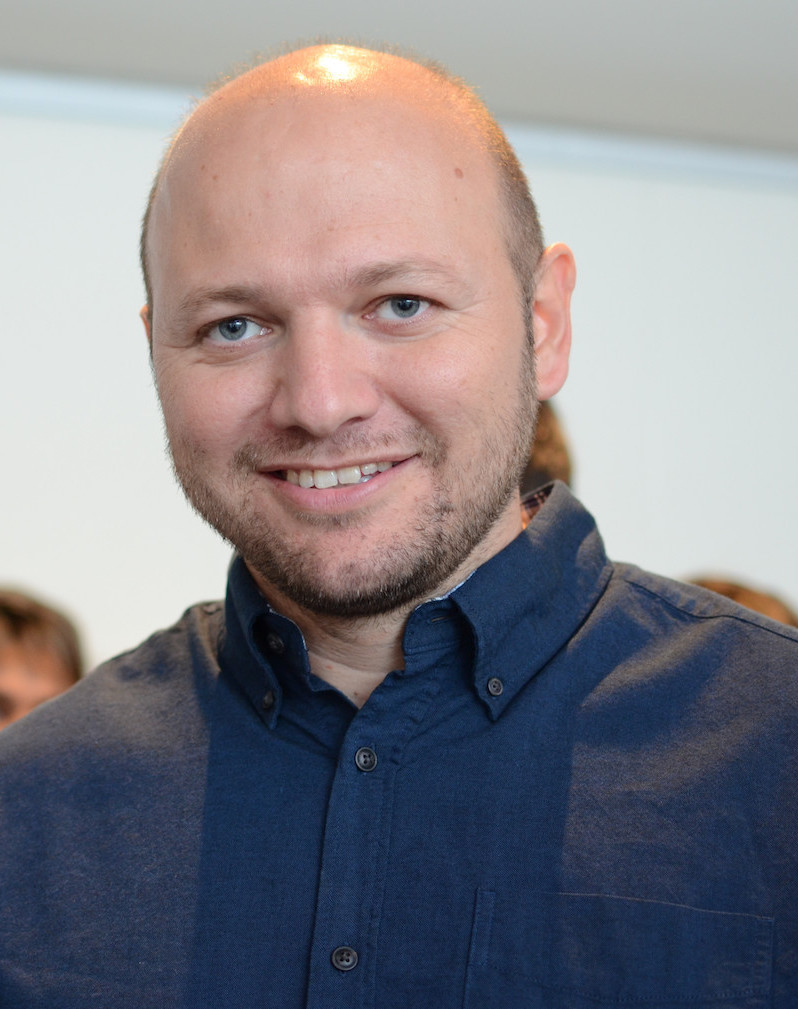}}]{Themistoklis Charalambous}
													received his BA and M.Eng in Electrical and Information Sciences from Trinity College, University of Cambridge. He completed his Ph.D. studies in the Control Laboratory, of the Engineering Department, University of Cambridge in 2010. Following his PhD, he worked as a research associate at Imperial College London, as a visiting lecturer at the Department of Electrical and Computer Engineering, University of Cyprus, and as a postdoctoral researcher at the Department of Automatic Control of the School of Electrical Engineering at the Royal Institute of Technology (KTH) and the Department of Electrical Engineering at Chalmers University of Technology. In January 2017, he joined the Department of Electrical Engineering and Automation, School of Electrical Engineering, Aalto University, as an Assistant Professor and in July 2020 he became an Associate Professor at Aalto University. Since September 2021, he is a tenure-track Assistant Professor at the Department of Electrical and Computer Engineering at the University of Cyprus and a Visiting Professor at Aalto University.
													
													His research involves distributed coordination and control, distributed decision making, and control of various resource allocation problems in complex and networked systems.
												\end{IEEEbiography}
												
											\end{document}